\documentclass[12pt]{amsart}

\usepackage[latin1]{inputenc}
\usepackage[T1]{fontenc}
\usepackage{amssymb,url}
\usepackage{amsthm}
\usepackage{amsfonts}
\usepackage[width=.75\textwidth,labelsep=period,font=footnotesize,labelfont={bf,sc},aboveskip=6pt]{caption}

\theoremstyle{plain}
\newtheorem{thm}{Definition}[section]
\newtheorem{theorem}[thm]{Theorem}
\newtheorem{lemma}[thm]{Lemma}

\newtheorem{proposition}[thm]{Proposition}
\newtheorem*{theorem*}{Theorem}

\theoremstyle{definition}
\newtheorem{definition}[thm]{Definition}

\newtheorem{example}[thm]{Example}

\theoremstyle{remark}
\newtheorem{remark}[thm]{Remark}
\newtheorem*{acknowledgements}{Acknowledgements}

\numberwithin{equation}{section}
\makeatletter
\newcommand{\low}{\@ifnextchar^{}{^{\vphantom x}}}
\newcommand{\high}{\@ifnextchar_{}{_{\vphantom I}}}
\makeatother

\makeatletter

\newcommand{\ltwoa}{\mathord{\mathpalette\tw@a<}}
\newcommand{\rtwoa}{\mathord{\mathpalette\tw@a>}}
\newcommand{\tw@a}[2]{\ooalign{\hfil$#1 #2$\hfil\crcr$#1 \Relbar\joinrel\Relbar$\crcr}}
\newcommand{\lthreea}{\mathord{\mathpalette\thr@a<}}
\newcommand{\rthreea}{\mathord{\mathpalette\thr@a>}}
\newcommand{\thr@a}[2]{\ooalign{\hfil$#1 #2$\hfil\crcr\raise.3ex\hbox{$#1 \relbar\mathrel{\mkern-4mu}\relbar$}\crcr
$#1 \relbar\mathrel{\mkern-4mu}\relbar$\crcr\lower.3ex\hbox{$#1 \relbar\mathrel{\mkern-4mu}\relbar$}\crcr}}
\newcommand{\Dyatop}[2]{\mathord{\mathpalette{\Dy@top{#1}{#2}}{}}}
\newcommand{\Dy@top}[3]{\overset{#3 #1}{#3 #2}}
\makeatother

\newcommand{\act}{\mathinner\cdot}
\DeclareMathOperator{\Ad}{Ad}
\DeclareMathOperator{\ad}{ad}
\def\al{\alpha}
\def\be{\beta}
\def\C{{\mathbb{C}}}
\def\Ct{{\C}^{\sn2}}
\newcommand{\cohom}[2]{\operatorname{cohom}\low_\sn{#1}{#2}}
\newcommand{\cross}[1]{\mathbin{{\times}\!}\low_\sn{#1}}
\def\de{\delta}

\DeclareMathOperator{\diag}{diag}
\DeclareMathOperator{\End}{End}
\def\G{\Lie{G}}
\def\ga{\gamma}
\newcommand{\go}[1]{\mathcal{#1}}
\def\gN{\go{N}}
\def\gZ{\go{Z}}
\def\geqs{\geqslant}
\def\H{{\mathbb{H}}}
\DeclareMathOperator{\Id}{Id}
\DeclareMathOperator{\im}{Im}
\newcommand{\inp}[3][]{\left\langle #2,#3\right\rangle\low_{#1}}
\newcommand{\inpsl}[2]{\left\langle #1,#2\right\rangle\low_{\sn{\Sl}}}
\def\ka{\kappa}
\def\la{\lambda}
\def\La{\Lambda}
\def\LA{\Lie{A}}
\def\LB{\Lie{B}}
\def\LC{\Lie{C}}

\newcommand{\DL}[1]{\lbrack\!\lbrack#1\rbrack\!\rbrack}
\def\LD{\Lie{D}}
\def\LE{\Lie{E}}
\def\leqs{\leqslant}
\def\LF{\Lie{F}}
\def\LG{\Lie{G}}
\def\LH{\Lie{H}}
\newcommand{\lie}[1]{\mathfrak{#1}}
\newcommand{\Lie}[1]{\textsl{#1}}
\def\lieg{\lie{g}}
\def\lieac{{\lie{a}}^{\sn\C}} 
\def\liegc{{\lie{g}}^{\sn\C}}
\def\Liegc{{\G}^{\sn\C}} 
\def\liehc{{\lie{h}}^{\sn\C}}
\def\liekc{{\lie{k}}^{\sn\C}}

\def\lierc{{\lie{r}}^{\sn\C}}
\newcommand{\List}[1]{{\left\{#1\right\}}}
\def\LK{\Lie{K}}
\newcommand{\Lm}[1]{\qq{L}{-}{#1}}
\newcommand{\Lp}[1]{\qq{L}{+}{#1}}
\def\lp{\lie{p}}
\def\LR{\Lie{R}}
\def\LS{\Lie{S}}
\def\LU{\Lie{U}}
\newcommand{\Lw}[1]{\ww{L}{#1}}
\newcommand{\lzero}{\text{\large 0}}
\def\O{\go{O}}
\newcommand{\ms}{{\raise1pt\hbox{$\scriptscriptstyle -$}}}
\def\N{{\mathbb{N}}}
\def\om{\omega}
\def\op{\oplus}
\def\ot{\otimes}
\newcommand{\Overline}[1]{
  \overline{\mathchoice{\displaystyle #1}{\textstyle #1}{\scriptstyle #1}{\scriptscriptstyle #1}}}
\def\pa{\partial}
\newcommand{\pms}{{\raise1pt\hbox{$\scriptscriptstyle\pm$}}}
\newcommand{\ps}{{\raise1pt\hbox{$\scriptscriptstyle +$}}}
\newcommand{\qq}[3]{#1^{\sn{#3}}\high_\sn{#2}}
\def\R{{\mathbb{R}}}
\DeclareMathOperator{\re}{Re}
\def\Rm{\R\low_\ms}
\def\Rp{\R\low_\ps}
\newcommand{\scr}{\scriptstyle}
\def\Se{{S}^{\sn1}}
\def\si{\sigma}
\def\Si{\Sigma}
\newcommand{\siEd}[1]{{\ww{E}{#1}}'}
\newcommand{\siXd}[1]{{\ww{X}{#1}}'}
\DeclareMathOperator{\Sl}{\lie{sl}}

\def\slc{{\lie{sl}(2,\C)}}
\def\slT{{{\slc}^{\sn3}}}

\newcommand{\sln}[1]{{\lie{sl}(#1,\C)}}
\newenvironment{smatrix}{\left(\smallmatrix}{\endsmallmatrix\right)}
\newcommand{\sn}[1]{{\raise1pt\hbox{$\scriptscriptstyle #1$}}}
\def\SO{\Lie{SO}}
\newcommand{\So}[1]{{\lie{so}(#1)}}
\newcommand{\son}[1]{{\lie{so}(#1,\C)}}
\def\SP{\Lie{Sp}}
\newcommand{\Span}[2][]{\operatorname{span}_{#1}{\{#2\}}}
\newcommand{\Sp}[1]{{\lie{sp}(#1)}}
\newcommand{\spn}[1]{{\lie{sp}(#1,\C)}}
\newcommand{\sscr}{\scriptscriptstyle}
\def\St{{S}^{\sn2}}
\newcommand{\su}[1]{{\lie{su}(#1)}}
\def\SU{\Lie{SU}}
\newcommand{\tablestrut}{\vrule height 11pt width 0pt}
\newcommand{\teta}{\tilde\eta}
\DeclareMathOperator{\tr}{Tr}
\def\TS{{\tablestrut}}
\def\vep{\varepsilon}
\def\we{\wedge}
\newcommand{\ww}[2]{#1\low_\sn{#2}}
\def\xia{{\xi\low_\sn A}}
\def\xib{{\xi\low_\sn B}}
\def\xiad{{{\xia}'}}
\def\xibd{{{\xib}'}}
\def\xiod{{{\ww{\xi}{1}}'}}
\def\xitd{{{\ww{\xi}{2}}'}}
\def\Z{{\mathbb{Z}}}
\def\ze{\zeta}

\begin{document}

\title[Cohomogeneity-Three HyperKähler Metrics]{Cohomogeneity-Three HyperKähler\\ Metrics on
  Nilpotent Orbits}

\author{Martin Villumsen}

\address{Department of Mathematics and Computer Science\\
  University of Southern Denmark\\
  Campusvej 55\\
  DK-5230 Odense M\\
  Denmark}

\email{mdv@imada.sdu.dk}


\begin{abstract}
  Let $\O$ be a nilpotent orbit in $\liegc$ where $\G$ is a compact, simple group and $\lieg=\Lie{Lie}(\G)$. It is known
  that $\O$ carries a unique $G$-invariant hyperK\"ahler metric admitting a hyperK\"ahler potential compatible with the
  Kirillov-Kostant-Souriau symplectic form. In this work, the hyperK\"ahler potential is explicitly calculated when $\O$
  is of cohomogeneity three under the action of $G$. It is found that such a structure lies on a one-parameter family of
  hyperK\"ahler metrics with $G$-invariant K\"ahler potentials if and only if $\lieg$ is $\Sp3$, $\su6$, $\So7$,
  $\So{12}$ or $\ww{\lie e}7$ and otherwise is the unique $\G$-invariant hyperK\"ahler metric with $\G$-invariant
  K\"ahler potential.
\end{abstract}


\maketitle
\section{Introduction}
Let $\G$ be a compact connected semisimple Lie group with Lie algebra $\lieg$ and let $\O$ be an adjoint nilpotent
$\Liegc$-orbit in the complexif\mbox{}ied Lie algebra $\liegc=\lieg\ot\C$. Then $\O$ carries a canonical
$\Liegc$-invariant complex symplectic structure known as the Kirillov-Kostant-Souriau symplectic structure (see
eg.~\cite{Besse:Einstein}). Using an identif\mbox{}ication of $\O$ with a set of solutions of Nahm's equations,
Kronheimer~\cite{Kronheimer} showed that $\O$ possesses a $\G$-invariant hyperK\"ahler metric.

Recall that a Riemannian manifold $(M,g)$ is \emph{hyperK\"ahler} if it admits three endomorphisms $I$, $J$ and $K$ of
the tangent bundle that satisfy the relations of the quaternions (i.e.~$I^\sn2=J^\sn2=-1$ and $IJ=K=-JI$), each
preserved by the Levi-Civita connection and each turning $g$ into a Hermitian metric. Necessarily, $I$, $J$, $K$ are
integrable so that $(M,g,I)$ etc.~are K\"ahler structures. As observed by Hitchin~\cite{Hitchin:monopoles}, replacing
the condition $\nabla I=\nabla J=\nabla K=0$ with $d\ww\om I=d\ww\om J=d\ww\om K=0$, the almost complex structures $I$,
$J$ and $K$ are still integrable; the two-forms $\ww\om I=g(I\cdot,\cdot)$ etc.~are called K\"ahler forms.  Note that if
one f\mbox{}ixes the complex structure $I$, then the non-degenerate and closed two-form $\om=\ww\om J+i\,\ww\om K$ is
holomorphic with respect to $I$, so $(M,I)$ carries a \emph{complex symplectic structure} def\mbox{}ined by $\om$.
Taking the $2n$'th power of $\om$ one gets a parallel section of the canonical bundle, revealing the well-known fact that
hyperK\"ahler metrics are Ricci-flat and Calabi-Yau.

By considering the cotangent bundle $T^{\sn *}\C P(n)$, Calabi~\cite{Calabi:Kaehler} gave the earliest non-trivial
examples of hyperK\"ahler metrics. Calabi's metric has symmetry group $\LU(n+1)$ and the space $T^{\sn *}\C P(n)$ is via
the complex moment map identif\mbox{}ied with an adjoint orbit in the semisimple Lie algebra $\lie{gl}(n+1,\C)$, and it
is of cohomogeneity one under the action of $\SU(n+1)$. Dancer \& Swann~\cite{DancerSwann:HK-metrics-of-cohom-one}
classif\mbox{}ied hyperK\"ahler metrics in cohomogeneity one with respect to a compact simple Lie group on manifolds of
dimension greater than four. There are two possibilities: the Calabi metric on $T^{\sn *}\C P(n)$ and nilpotent orbits
of cohomogeneity one in complex simple Lie algebras. A feature of all of these metrics is that they admit an invariant
K\"ahler potential. This aspect makes it natural to consider nilpotent orbits of a given cohomogeneity and with
invariant hyperK\"ahler metrics that admit K\"ahler potentials.

A hyperK\"ahler potential is a function that is simultaneously a K\"ahler potential of each complex structure compatible
with the hyperK\"ahler metric. By Swann~\cite{Swann:HK-quaternionicKahler} it is known that the hyperK\"ahler structure
of Kronheimer on a nilpotent orbit admits a $\G$-invariant hyperK\"ahler potential. In fact, it is the unique
$\G$-invariant hyperK\"ahler potential on $\O$, as shown by Brylinski~\cite{Brylinski}.

In this paper we study nilpotent orbits $\O$ of cohomogeneity three under the action of the compact group $G$. We
restrict ourselves to the case when $\lieg$ is simple and consider $G$-invariant hyperK\"ahler structures on $\O$
compatible with the Kirillov-Kostant-Souriau symplectic structure, that admit a $\G$-invariant K\"ahler potential.
Expressing the metric in terms of the Killing form and the Lie bracket we explicitly f\mbox{}ind the unique
$\G$-invariant hyperK\"ahler potential on each orbit. We show that this potential lies in a larger (one-dimensional)
family of $\G$-invariant hyperK\"ahler metrics with $\G$-invariant K\"ahler potential if and only if $\lieg$ is $\Sp3$,
$\su6$, $\So7$, $\So{12}$ or $\ww{\lie e}7$.

For orbits of cohomogeneity two Kobak \& Swann~\cite{KobakSwann:HK-pot-in-cohom-two} found similar results.

It is essential that when $\G$ is not $\SO(7)$ each element of a nilpotent orbit of cohomogeneity three can be moved
into a subalgebra of three commuting $\si$-invariant copies of $\slc$, where $\si$ is the compact real structure. This
feature is shown using an interesting relation with rank three Hermitian symmetric spaces and it allows one to describe
the geometry using representations of $\slc$. On the other hand, the cohomogeneity three nilpotent orbit in
$\O_{(3,2^\sn2)}\subset\son7$ is special: for generic $X\in\O_{(3,2^\sn2)}$ the algebra generated by $X$ and $\si X$ is
all of $\son7$, so a separate direct approach has to be used. We f\mbox{}ind that there is only a one-dimensional family
of hyperK\"ahler metrics with K\"ahler potential. This is surprising for the following reason. It is known that
$\O_{(3,2^\sn2)}\subset\son7$ has a one-to-one correspondence with the cohomogeneity-two nilpotent orbit in
$\O_{(2^\sn4)}\subset\son8$ that carries a one-parameter family of $\SO(8)$-invariant hyperK\"ahler metrics with
$\SO(8)$-invariant K\"ahler potential. Consequently, reducing the symmetry group from $\SO(8)$ to $\SO(7)$ give us no
extra families even though the cohomogeneity changes from two to three.

\section{Geometry of Nilpotent Orbits}
\label{Section:geo-nilpotent-orbit}
In this section we review some known facts on nilpotent orbits and establish basic notation. We introduce the techniques
used to compute cohomogeneities on nilpotent orbits and hopefully clarify the results given in
Tables~\ref{Table:Big-table-classical} and~\ref{Table:Big-table-execptional} via examples.

\subsection{Principal Orbits and Cohomogeneities}
If $\G$ is a Lie group acting on some space $M$ then two orbits $\G\act x,\G\act y\subseteq M$ are said to have the same
\emph{orbit~type} if the isotropy subgroups $\ww\G x$ and $\ww\G y$ are conjugate in $\G$. If $M$ is a manifold the
following fundamental theorem shows the existence of a maximal orbit type.

\begin{theorem}[see~\cite{Bourbaki1982,Bredon}]\label{theorem:principal-orbits}
  Let $M$ be a manifold admitting a dif\mbox{}ferentiable action of a compact Lie group $\G$. Suppose the orbit space
  $M/\G$ is connected. Then there exists a maximal orbit type $(\LH)$. That is, $\LH$ is conjugate to a subgroup of each
  isotropy group. The union $\ww M{(\LH)}$ of the orbits of type $(\LH)$ is open and dense in $M$ and $\ww M{(\LH)}/\G$
  is connected. \qed
\end{theorem}

An orbit of type $(\LH)$, with $\LH$ as in Theorem~\ref{theorem:principal-orbits}, is called a \emph{principal~orbit}
and the \emph{cohomogeneity} of the action of $\G$ on $M$ is def\mbox{}ined to be the codimension of a principal
orbit.

With $\G$ compact each point $p\in M$ has a slice. This is a submanifold $S$ such that $\G S\subset M$ is open and there
exists a $\G$-equivariant isomorphism of the $\G$-spaces $\G S$ and $\G\cross{\LH} S$, where $H$ is the isotropy group
$\ww\G p$. Note in particular that if $M$ is a vector space and $\G$ a compact Lie group then $\G S\cong\G\cross{\LH}W$
(cf.~\cite{Bredon}), where $W$ is the normal space of the tangent space $\ww Tp(\G\act p)$. This is of great importance
to us because we thereby get the equality
\begin{equation}
  \label{eq:cohom-G=gohom-H}
  \cohom\G{\G S}=\cohom{\LH}W\text{.}
\end{equation}

\begin{example}
  \label{ex:cohom-calc-Bredon}
  Let us use~\eqref{eq:cohom-G=gohom-H} to compute the cohomogeneity of some group representations. We start of\mbox{}f
  with a simple example.
  
  \textup{(i)} Let $V^\sn n$ denote the $\LU(n)$-module $\C^\sn n$ via the standard action and let $\List{\ww ej}$ be the
  standard basis of $V^\sn n$.  The stabilizer in $\LU(n)$ of $\ww e1$ is
  \begin{equation*}
    \LH=\ww{\LU(n)}{\ww e1} =\List{\begin{smatrix}1&0\\0&A\end{smatrix}|A\in\LU(n-1)}\cong\LU(n-1)\text{.}
  \end{equation*}
  As an $\LH$-module we now have $\ww T{\ww e1}V^\sn n=V^\sn n=\C\op V^\sn{n-1}$, and in turn $\lie u (n)\cong \R\op
  V^\sn{n-1}\op\lie u(n-1)$, where $\R$, $\C$ are trivial modules. Thus, $ \ww T{\ww e1}\left(\LU(n)\cdot \ww
    e1\right)\cong\R\op V^\sn{n-1}$ and we conclude that
  \begin{equation*}
    \cohom{\LU(n)}{V^\sn n}=\cohom{\LH}{\R}=1\text{.}
  \end{equation*}
    
  Since the unitary group acts transitively on the sphere, the above result should come as no surprise. F\mbox{}inding
  the cohomogeneity of the following representations is less trivial.
  
  \textup{(ii)} Consider the $U(n)$ module $\La^\sn2 V^\sn n=\Span[\sn\C]{\ww ei\we \ww ej|1\leqslant i<j\leqslant
    n}$. The stabilizer of $\ww e1\we\ww e2$ is
  \begin{equation*}
    \begin{split}
      \LH=\ww{\LU(n)}{\ww e1\we\ww e2} &=\List{\begin{smatrix}A&0\\0&B\end{smatrix}|A\in\SU(2),B\in\LU(n-2)}\\
      &\cong\SU(2)\times\LU(n-2)\text{.}
    \end{split}
  \end{equation*}
  As an $\LH$-module we have $V^\sn n=\Se V^\sn2\op V^\sn{n-2}$, giving us
  \begin{equation*}
    \ww T{\ww e1\we\ww e2}\La^\sn2 V^\sn n=\La^\sn2 V^\sn n=\La^\sn2\left(\Se V^\sn2\right) \op \left(\Se V^\sn2\ot
    V^\sn{n-2}\right)\op\La^\sn2 V^\sn{n-2}\text{.}
  \end{equation*}
  Moreover $\lie u(n)=\lie u(2)\op\lie u(n-2)\op\left(\Se V^\sn2\ot V^\sn{n-2}\right)$, so that
  \begin{equation*}
   \ww T{\ww e1\we\ww
    e2}\left(\LU(n)\cdot(\ww e1\we\ww e2)\right)\cong\R\op(\Se V^\sn2\ot V^\sn{n-2})\text{.}
  \end{equation*}
  The cohomogeneity is now
  \begin{equation*}
   \begin{split}
      \cohom{\LU(n)}{\La^\sn2 V^\sn n}&=\cohom{\LH}{\left(\La^\sn2\R^\sn2\op\La^\sn2 V^\sn{n-2}\right)}\\
      &=1+\cohom{\LU(n-2)}{\La^\sn2}{V^\sn{n-2}}\\
      &= \left\{
        \begin{array}{ll}
          n/2     &\  n \ \mbox{\text{even}}\\
          (n-1)/2 &\  n \ \mbox{\text{odd}}
        \end{array}\right.\text{,}
    \end{split}
  \end{equation*}
 by induction.
 
  \textup{(iii)} Consider the symmetric product $\St V^\sn3$. The stabilizer in $\LU(3)$ of $\ww e1\vee\ww e1$ is
  \begin{equation*}
    \LH=\List{
      \begin{smatrix}
        \pm 1&0\\0&A
      \end{smatrix}
      |A\in\LU(2)}\cong\ww\Z2\times\LU(2)\text{.}
  \end{equation*}
  As $\Lie{H}$-modules we have $V^\sn3=\C\op V^\sn2$, $\St V^\sn3\cong\C\op V^\sn2\op \St V^\sn2$ and $\lie
  u(3)=\R\op V^\sn2\op\lie u(2)$, so that $\ww T{\ww e1\vee\ww e1}\left(U(3)\cdot \ww e1\vee\ww e1\right)\cong\R\op
  V^\sn2$.  The cohomogeneity is thus
  \begin{equation*}
    \begin{split}
      \cohom{\LU(3)}{\St V^\sn3}&=\cohom{\LH}{\left(\R\op \St V^\sn2\right)}\\
      &=1+\cohom{\LU(2)}{\St V^\sn2}=3\text{.}
    \end{split}
  \end{equation*}
\end{example}

\subsection{Adjoint Nilpotent Orbits}
Letting $\si$ denote the conjugation map of $\liegc$ with respect to the compact real form $\lieg$, the following
def\mbox{}ines a positive def\mbox{}inite symmetric bilinear form on the real Lie algebra $(\liegc)^\sn\R$
\begin{equation*}
  (X,Y)\mapsto \inp{X}{\si Y},\quad X,Y\in(\liegc)^\sn\R\text{,}
\end{equation*}
where $\inp\cdot\cdot$ is the negative of the Killing form.

To ease notation we recursively def\mbox{}ine a multiple bracket $\DL{\dots}$ on the Lie algebra $(\liegc,[\cdot,\cdot])$,
\begin{equation*}
  \begin{split}
  \DL{\ww X1\ww X2}&=[\ww X1,\ww X2],\\
  \DL{\ww X1\ww X2\dots\ww Xk}&=[\ww X1,\DL{\ww X2\dots\ww Xk}],\quad\ww Xj\in\liegc\text{.}
  \end{split}
\end{equation*}

Via the adjoint action each element $A\in\liegc$ generates a vector f\mbox{}ield $\xia$ on a nilpotent orbit $\O$.  This
is given by $\xia\ww|X=\DL{AX}$. When there is no danger of confusion, the subscript $X$ will be suppressed throughout
this paper. An important property of these vector f\mbox{}ields is the identity $\DL{{\xia}{\xib}}=-\ww\xi{\DL{AB}}$.
As the tangent space of $\O$ at $X$ is $\List{\DL{AX}| A \in \liegc}$, we see that the complex structure $I$ is
characterized by $I\xia=\ww\xi{i A}$.

The nilpotent orbit $\O\subset\liegc$ is a complex submanifold with complex structure $I$ which is inherited from the
natural embedding. On $\O$ the \emph{Kirillov-Kostant-Souriau} symplectic form is given by,
\begin{equation*}
  \ww{\Si(\xia,\xib)}{X}=\inp{\DL{AB}}{X}\text{.}
\end{equation*}
We say that a hyperK\"ahler structure $(g,I,J,K)$ is compatible with $(\O,\Si,I)$ if it satisf\mbox{}ies the equation
$\Si=\ww\om J+i\,\ww\om K$.

The compact group acts naturally on the non-compact orbit $\O$ and we shall have a particular interest in
the situation with three parameters transverse to the action, that is when $\O$ is a nilpotent orbit of cohomogeneity
three. However, we notice that
\begin{remark}
  \label{remark:R-times-G}
  If $\cohom\G\O=1$ then $\Rp\times\G$ acts transitively on $\O$.
\end{remark}

Orbits in classical Lie algebras are essentially determined by Jordan normal forms of elements in the standard
representation, so they are specif\mbox{}ied by partitions. This characterization is not useful for orbits in
exceptional Lie algebras, which are best described by a weighted Dynkin diagram (see~\cite{CollingwoodMcGovern}).

\subsection{Standard Triples}
\label{Section:StandardTriple}
A set $\List{H,X,Y}$ of nonzero elements in $\liegc$ is called a \emph{standard triple} if the following bracket
relations hold,
\begin{equation*}
  \DL{HX}=2\,X, \quad\DL{HY}=-2\,Y\quad\text{and}\quad\DL{XY}=H\text{.}
\end{equation*}  
Since the element $H$ is a semisimple element of the Lie subalgebra spanned by $\List{H,X,Y}$, which in turn is
isomorphic to the simple Lie algebra $\slc$, it must be semisimple as an element of $\liegc$.  Similarly $X$, $Y$ are
nilpotent elements of $\liegc$.

The Jacobson-Morozov theorem is of fundamental importance here because it provides access to the representation theory
of $\slc$; a detailed proof may be found in~\cite{CollingwoodMcGovern}.
\begin{theorem*}[Jacobson--Morozov]
  Let $\liegc$ be a semisimple Lie algebra. If $X$ is a nonzero nilpotent element of $\liegc$, then $X$ embeds into a
  standard triple $\List{H,X,Y}$ of $\liegc$. \qed
\end{theorem*}

Given a compact real form one may obtain a $\si$-invariant $\slc$-subalgebra associated to a nilpotent orbit.

\begin{proposition}[Borel]\label{Borel}
  Let $\O$ be a nilpotent orbit in $\liegc$ and let $\si$ be the conjugation of a compact real form. Then $\O$ contains
  an element X such that $\List{\DL{X(\si X)},X,\si X}$ spans a Lie subalgebra isomorphic to $\slc$.
\end{proposition}

\begin{proof}
  See~\cite{KobakSwann:HK-pot-in-cohom-two}.
\end{proof}

\subsection{The Closure Ordering}
There exists a partial order on the set of nilpotent orbits in $\liegc$ def\mbox{}ined by $\O'\preceq\O$ if and only if
$\O'\subset\Overline\O$. Notice that the boundary of any orbit is the union of orbits of smaller dimensions,
\begin{equation}
  \label{eq:Boundary-relation}
  \Overline{\O'}=\ww{\bigcup}{\O\preceq\O'}\O\text{.}
\end{equation}

On the level of cohomogeneity it is a pleasure to note the following Theorem by Dancer and
Swann~\cite{DancerSwann:qK-cohom-one}.

\begin{theorem}
  \label{the:DancerSwann}
  Let $\G$ be a compact simple Lie group. If $\ww\O1$ and $\ww\O2$ are nilpotent orbits with $\ww\O1\precneqq\ww\O2$,
  then
  \begin{displaymath}
    \cohom\G{\ww\O1}<\cohom\G{\ww\O2}\text{.}
  \end{displaymath}
\qed
\end{theorem}

\begin{example}
  Any element $X$ of the nilpotent orbit in $\sln n$ characterized by the Jordan normal form $d=(3,1^\sn{n-3})$ with
  $n\geqs4$ is $\SU(n)$-conjugated to an element of the form
\begin{equation*}
  \begin{smatrix}
    0&s&\la&r& \\
    0&0&t&0&\\
    0&0&0&0&\\
    0&0&0&0&\\
    &&&&\text{\large 0}\\
  \end{smatrix}
  \quad r,s,t\in\Rp\quad\la\in\C\text{.}
\end{equation*}
Consequently the cohomogeneity of $\ww\O d\subset\sln n$ under the action of the compact group $\SU(n)$ is at most $5$.
Whenever $n\geqslant4$, $\ww\O d$ is neither minimal nor next-to-minimal (see Section~\ref{sec:minimal-orbit}), so by
Theorem~\ref{the:DancerSwann} the cohomogeneity is $3$,$4$ or $5$. Evidently, for $n=4$, the cohomogeneity is no less
than $5$. But for $n\geqs5$ we need more information (see Example~\ref{ex:comom-calc}).
\end{example}

\subsection{The Minimal Orbit}
\label{sec:minimal-orbit}
When $\liegc$ is simple there is a unique nonzero minimal orbit $\ww\O{\text{min}}\neq\List0$ with respect to the
closure ordering (see~\cite{CollingwoodMcGovern}). In fact $\ww\O{\text{min}}$ corresponds to a highest root, so it is
the orbit through a nonzero element in a highest root space.

We label the orbits above the minimal orbit via
\begin{definition}
  \textup{(i)} The \emph{minimal nilpotent orbit} $\ww\O{\text{min}}$ is the unique nilpotent orbit characterized by
  $\ww\O{\text{min}} \preceq \O$ for any nonzero orbit $\O$. The minimal orbit will also be referred to as \emph{the
    height one} orbit.
  
  \textup{(ii)} The \emph{height} of a nilpotent orbit $\O\succneqq\ww\O{\text{min}}$ is def\mbox{}ined by induction.
  Let $\ww S\O$ be the set of nilpotent orbits $\O'\precneqq\O$ for which there is no nilpotent orbit $\O''$ satisfying
  $\O'\precneqq\O''\precneqq\O$. Let $\ww\O0\in\ww S\O$ be an orbit of minimal height $n$. Then the height of $\O$ is
  $n+1$.
\end{definition}
\noindent
A a nilpotent orbit of height two will also be referred to as \emph{next-to-minimal}. Notice that the Lie algebra
may possess more than one next-to-minimal orbit, as well as orbits of height three.

The minimal orbit is in fact the unique nilpotent orbit of cohomogeneity one under the action of the compact group
$\G$. If $\al^{\sn\#}$ is a highest root we get a \emph{highest root decomposition} by

\begin{proposition}
  \label{prop:DancerSwann-cohom1}
  \textup{(i)} As an $\ww{\slc}{\al^\sn{\#}}$-module, the Lie algebra $\liegc$ decomposes as
  \begin{equation}
    \label{eq:wolf-lie-algebra-level}
    \liegc \cong \ww{\slc}{\al^\sn{\#}} \op \liekc \op (\Se \ot V)\text{,}
  \end{equation}
  where $\liekc$ is the centralizer of $\ww{\slc}{\al^\sn{\#}}$, $\Se=\Se\Ct$ is the fundamental
  $\ww{\slc}{\al^\sn\#}$ representation, and $V$ is a trivial $\ww{\slc}{\al^\sn{\#}}$-module and a non-trivial
  $\liekc$-module.

  \textup{(ii)}Under the action of the compact group $\G$ the nilpotent orbit $\ww\O{\text{min}}$ is of cohomogeneity
  one.
\end{proposition}

\begin{proof}
  See Kobak and Swann~\cite{KobakSwann:HK-geo-wolf-spaces}.
\end{proof}

The compact homogeneous space $W(\G)=\G/(\SU(2)\LK)$ corresponding to a highest root decomposition, is a symmetric
space with a quaternionic K\"ahler structure, cf.~\cite{Swann:HK-quaternionicKahler}. Such manifolds are known as
\emph{Wolf spaces} and are listed in Besse~\cite[p.~409]{Besse:Einstein}.
\begin{remark}
  \label{remark:real-dim-of-wolf-space}
  Via~\eqref{eq:wolf-lie-algebra-level} the complexif\mbox{}ication of the tangent space of $W(\G)$ is isomorphic to
  $\Se\ot V$, so
  \begin{equation}
    \label{eq:dim-of-V}
    \ww{\dim}{\C}V=\tfrac12\ww{\dim}{\R}W(\G)\text{.}
  \end{equation}
\end{remark}

\subsection{The Beauville Bundle}
\label{sec:Beauville}
A useful approach to the computation of the cohomogeneity of a nilpotent orbit $\O$ is to consider the Beauville bundle
$\go N(\O)$~\cite{Beauville:fano}. $\gN(\O)$ is a vector bundle that contains an open orbit which is identif\mbox{}ied
with $\O$.  Following Kobak \& Swann~\cite{KobakSwann:HK-pot-in-cohom-two} we shall outline a procedure for
f\mbox{}inding cohomogeneities.

Let $\O$ be a nilpotent orbit in the semisimple Lie algebra $\liegc$. Let $\si$ be the conjugation with respect to the
compact real form $\lieg$, and choose a standard triple $\List{H,X,Y=-\si X}$ of a real $\slc$ subalgebra with $X\in\O$.
By the representation theory of $\slc$ the eigenvalues of $\ww\ad H$ are integers. For each $j\in\N$ we let $\liegc(j)$
be the corresponding eigenspace of $\ww\ad H$ and set $\lp=\ww\bigoplus{i\geqslant 0} \liegc(i)$ and $\lie n =
\ww\bigoplus{i\geqslant 2} \liegc(i)$. Then $\lp$ is a parabolic subalgebra, $\lie n$ is a nilpotent subalgebra and
$X\in\lie n$. Let $\Lie P$ be the corresponding parabolic subgroup $\Lie P=\List{g\in\Liegc|\ww\Ad g(\lp)=\lp}$.  This
is a closed complex Lie subgroup of $\Liegc$ with Lie algebra equal to $\lp$. In particular $\Liegc/\Lie P$ is a
homogeneous space. Clearly $\lie n$ is an ideal of $\lp$ and there is a natural action of $\Lie P$ on $\Liegc\times\lie
n$,
\begin{equation*}
  p\cdot (g,N)\mapsto(g p^\sn{-1},\ww\Ad p N),\quad p\in\Lie P,\:g\in\Liegc,\:N\in\lie n\text{.}
\end{equation*}
The \emph{Beauville bundle} $\gN(\O)$ is def\mbox{}ined to be the vector bundle
\begin{equation*}
  \gN(\O) =\Liegc\cross{\Lie P} \lie n
\end{equation*}
over $\Liegc/\Lie P$ with f\mbox{}ibre $\lie n$. Letting $\pi$ denote the quotient map $\Liegc\times\lie n\to\gN(\O)$
we have the following equality of stabilizers, $\qq\G{\pi(e,N)}\C=\ww{\Lie P}N$, $N\in\lie n$. If $\tilde{\O}$ denotes
the $\Liegc$-orbit in $\gN(\O)$ then
\begin{equation*}
  \O\cong\Liegc/\qq{\G}{X}{\C}=\Liegc/\qq{\G}{{\pi(e,X)}}{\C}\cong\tilde\O\text{.}
\end{equation*}
That is, we have a $\Liegc$-equivariant isomorphism of the orbits $\O$ and $\tilde\O$, def\mbox{}ined by $g\cdot X
\mapsto g\cdot\pi(e,X)$, $g\in\Liegc$.

Let $\lie k\subseteq\lieg$ be the subalgebra satisfying $\lie k^\sn\C=\liegc(0)$ and let $\LK\subseteq\G$ be the
connected Lie subgroup with Lie algebra $\lie k$. Now $\Liegc/P=\G/\LK$ so $\gN(\O)=\G\cross{\LK}\lie n$. By
Theorem~\ref{theorem:principal-orbits} the union of the principal orbits in the $\G$-space $\gN(\O)$ forms a dense set.
Since $\tilde{\O}$ is an open subset, a principal orbit of the $\G$-space $\tilde\O$ must be a principal orbit of the
$\G$-space $\gN(\O)$. Using the $\G$-equivariant isomorphism we conclude
\begin{equation}
  \label{eq:cohom-formula}
  \cohom\G\O=\cohom\G{\tilde\O}=\cohom\G{\gN(\O)}=\cohom\LK{\lie n}\text{.}
\end{equation}
To f\mbox{}ind the cohomogeneity of $\lie n$ under the action of $\LK$, one merely has to use equation
\eqref{eq:cohom-G=gohom-H}.

\begin{example}
  \label{ex:comom-calc}
  We illustrate the \emph{Beauville bundle method} by considering the nilpotent orbit $\O$ in $\sln n$ with Jordan
  normal form $(3,1^\sn{n-3})$. For $n\geqs4$ this is an orbit of height three. We take the obvious choice of a standard
  triple (see Section~\ref{sec:explicit_standard_triples}) and obtain a $K=\ww{\LU(1)}+\times\ww{\Lie
    U(1)}-\times\SU(n-2)$ representation
  \begin{equation*}
    \lie n=\Lp2\Lm3\op\Lw+\Lm4\Se\C^\sn{n-2} \op \Lw+\Lm{-1} \Se\C^\sn{n-2}\text{.}
  \end{equation*}
  F\mbox{}ixing the f\mbox{}irst circle action, the problem is reduced to f\mbox{}inding the cohomogeneity of the
  $\ww{\LK}2:=\LU(1)\SU(n-2)$ representation
  \begin{equation*}
    \ww{\lie n}2:= L^\sn4 \Se\C^\sn{n-2} \op L^\sn{-1} \Se\C^\sn{n-2}\text{,}
  \end{equation*}
  Using the method described in Example~\ref{ex:cohom-calc-Bredon} one f\mbox{}inds that
  \begin{equation*}
    \begin{split}
      \cohom{\SU(n)}\O&=1+\cohom{\ww{\LK}2}{\ww{\lie n}2}=3 + \cohom{\SU(n-3)}\Se\C^\sn{n-3}\\
      &=\left\{
        \begin{array}{ll}
          5 &\ n=4\\
          4 &\ n\geqs5
        \end{array}
      \right.\text{.}
    \end{split}  
  \end{equation*}
  For $n=3$, the orbit $(3)$ is the regular orbit of $\sln3$. This is next-to-minimal and of cohomogeneity four
  (see~\cite{DancerSwann:qK-cohom-one}).
\end{example}

\section{Orbits of Height Three}
By work of Dancer \& Swann~\cite{DancerSwann:qK-cohom-one} all next-to-minimal orbits but one are of cohomogeneity two;
the next-to-minimal orbit of $\sln3$ (the regular orbit) is of cohomogeneity four. In the view of
Theorem~\ref{the:DancerSwann}, the nilpotent orbits of cohomogeneity three, if any, must therefore be of height three in
the partial order. So indeed the orbits of Table~\ref{tab:min-next-thrid} are our candidates in the search for orbits of
cohomogeneity three (recall that $\ww\LA3=\ww\LD3$, $\ww\LB2=\ww\LC2$).
\begin{table}[h]
  \begin{center}
    \begin{tabular}[t]{lc}
      \hline\hline
      \TS Type           & Height three orbit       \\
      \hline\hline
      \TS $\ww\LA n$     &                          \\
      \TS $n\geqslant 4$ & $(3,1^\sn{n-2})$         \\
      \TS $n\geqslant 5$ & $(2^\sn3,1^\sn{n-5})$    \\
      \hline
      \TS $\ww\LB n$     &                          \\
      \TS $n=2$          & $(5)$                    \\
      \TS $n\geqslant 3$ & $(3,2^\sn2,1^\sn{2n-6})$ \\
      \TS $n\geqslant 6$ & $(2^\sn6,1^\sn{2n-11})$  \\
      \hline
      \TS $\ww\LC n$     &                          \\
      \TS $n\geqslant 3$ & $(2^\sn3,1^\sn{2n-12})$  \\
      \hline
      \TS $\ww\LD n$     &                          \\
      \TS $n=3$          & $(3^\sn2)$               \\
      \TS $n\geqslant 4$ & $(3,2^\sn2,1^\sn{2n-7})$ \\
      \TS $n\geqslant 6$ & $(2^\sn6,1^\sn{2n-12})$  \\
      \hline\hline
    \end{tabular}
    \qquad
    \begin{tabular}[t]{lc}
      \hline\hline
      \TS Type      & Height three orbit     \\
      \hline\hline
      \TS $\ww\LG2$ & $\scr 2\rthreea0$      \\
      \TS $\ww\LF4$ & $\scr 01\rtwoa00$      \\
      \TS $\ww\LE6$ & $\scr 00\Dyatop0100$   \\
      \TS $\ww\LE7$ & $\scr 200\Dyatop0000$ , $\scr 000\Dyatop0010$\\
      \TS $\ww\LE8$ & $\scr 0100\Dyatop0000$ \\
      \hline\hline
    \end{tabular}
  \end{center}
  \caption[Orbits of height three in simple Lie algebras]{Orbits of height three in simple Lie algebras.}
  \label{tab:min-next-thrid}
\end{table}

\subsection{Explicit Standard Triples}
\label{sec:explicit_standard_triples}
The way of f\mbox{}inding the cohomogeneity via the methods described in~\S\ref{sec:Beauville} is based upon the
Jacobson-Morozov Theorem. Proposition~\ref{Borel} ensures the existence of a real $\si$-invariant $\slc$
subalgebra associated to the nilpotent orbit $\O$. In this section we shall explicitly list the elements $X$ of a standard
triple $\List{\DL{XY},X,Y}$ with $\si X=-Y$, associated to nilpotent orbits of height three in classical Lie algebra.

Note that if $\lieg$ is an exceptional Lie algebra, the weighted Dynkin diagram provides the necessary information.
  
The complex simple Lie algebras of type $\ww\LB{(n-1)/2}$ and $\ww\LD{n/2}$ can be realized as the set $\son
n= \List{Z\in\lie{gl}(n,\C)|ZB+BZ^\sn t=0}$ where $B$ is a symmetric invertible matrix. The standard compact real form
of $\lie{so(n,\C)}$ is $\si\colon Z\mapsto B\Overline ZB=-{\Overline{Z^\sn t}}$. The standard choice for $B$ is the
identity, but in this context we choose $B$ dif\mbox{}ferently. For orbits with Jordan normal form
$(2^\sn6,1^\sn{n-12})$ we take
\begin{equation}
  \label{eq:B12-Ik}
  B=
  \begin{pmatrix} 
    \ww B{12}&\lzero \\
    \lzero&\ww I{n-12}
  \end{pmatrix}
\end{equation}
where $\ww{B}{12}$ is the $(12\times12)$ matrix with $1$'s down the anti-diagonal and $0$'s elsewhere. In this picture
we let
\begin{equation}
  \label{eq:X-222222-0}
  X=\diag\Bigl(\begin{smatrix} 0&1\\0&0 \end{smatrix},\begin{smatrix} 0&1\\0&0 \end{smatrix},
  \begin{smatrix} 0&1\\0&0 \end{smatrix}, \begin{smatrix} 0&-1\\0&0 \end{smatrix},
  \begin{smatrix} 0&-1\\0&0 \end{smatrix},\begin{smatrix} 0&-1\\0&0 \end{smatrix},\ww\lzero{n-12} \Bigr)\text{.}
\end{equation}
For the orbit in $\ww\LB2$ with Jordan normal form $(5)$, $B$ is the $(5\times5)$ matrix with $1$'s down the
anti-diagonal, and our choice of $X$ is
\begin{equation*}
  X=
  \begin{smatrix}
    0&\sqrt2(1-i)&0&0&0\\
    0&0&-\sqrt6i&0&0\\
    0&0&0&\sqrt6i&0\\
    0&0&0&0&\sqrt2(i-1)\\
    0&0&0&0&0
  \end{smatrix}\text{.}
\end{equation*}
For the orbit with Jordan normal form $(3^\sn2)$ in $\ww\LD3$, $B$ is the $(6\times6)$ matrix with with $1$'s down
the anti-diagonal, and
\begin{equation*}
  X=i\sqrt2
  \begin{smatrix}
    0&1&0&&&\\
    0&0&1&&&\\
    0&0&0&&&\\
    &&&0&-1&0\\
    &&&0&0&-1\\
    &&&0&0&0
  \end{smatrix}\text{.}
\end{equation*}
When considering the orbit in $\ww\LB3$ characterized by the Jordan type $(3,2^\sn2)$ we choose
\begin{equation*}
  B=
  \begin{pmatrix}
    1&0&0\\
    0&0&\ww I3\\
    0&\ww I3&0
  \end{pmatrix} \quad\text{and}\quad
  X=
  \begin{pmatrix}
    0&0&v^\sn t\\
    -v&0&A\\
    0&0&0
  \end{pmatrix}\text{,}
\end{equation*}
where $v=(\sqrt2,0,0)$ and $A=\begin{smatrix} 0&0&0 \\ 0&0&1 \\ 0&-1&0 \end{smatrix}$. By adding trivial blocks, as
in~\eqref{eq:B12-Ik} and~\eqref{eq:X-222222-0}, this is easily generalized to the orbits $(3,2^\sn2,1^\sn{n-7})$ in
$\ww\LB{(n-1)/2}$, $\ww\LD{n/2}$.

For type $\ww\LA n=\sln{n+1}$ there are two orbits to consider, characterized by the partitions
$(2^\sn3,1^\sn{n-5})$ and $(3,1^\sn{n-2})$. Suitable choices of $X$ are evident,
\begin{equation*}
  X=
  \begin{smatrix}
    0&1&&&&&\\
    0&0&&&&&\\
    &&0&1&&&\\
    &&0&0&&&\\
    &&&&0&1&\\
    &&&&0&0&\\
    &&&&&&\lzero
  \end{smatrix}
  \quad\text{and}\quad
  X=
  \begin{smatrix}
    0&1&0&\\
    0&0&1&\\
    0&0&0&\\
    &&&\lzero
  \end{smatrix}\text{,}
\end{equation*}
respectively. These have the desired properties with respect to the standard compact real form $\si\colon Z\mapsto
-\Overline{Z^\sn t}$.

For type $\ww\LC n=\spn n$ there is one orbit to consider. We take
\begin{equation*}
  X=
  \begin{pmatrix}
    0&A\\
    0&0
  \end{pmatrix},
  \quad\text{with}\quad
  A=
  \begin{smatrix}
    1&0&0&\\
    0&1&0&\\
    0&0&1&\\
    &&&\lzero
  \end{smatrix}\text{.}
\end{equation*}

\subsection{The Cohomogeneities}
With the choices of the previous section we are able to f\mbox{}ind the centralizer $\lie t$ and
representation $\lie n$ of~\S\ref{sec:Beauville} in each case. This leads to Tables~\ref{Table:Big-table-classical}
\&~\ref{Table:Big-table-execptional} which exhaust the complete list of adjoint nilpotent orbits of cohomogeneity three in
simple Lie algebras.

\begin{table}[tb]
  {\small
      \begin{tabular}{lll}
        \hline\hline
        \TS Type$\backslash$Orbit    & $\lie k$$\backslash$$\lie n$ &  \makebox[0.5cm][r]{$\cohom\G\O$}\\
        \hline\hline
        \TS$\ww\LA5$                 & $\R+\ww{\su3}a+\ww{\su3}b$                  & $3$ \\
        \TS$(2^\sn3)$                & $L\End(\C^\sn3)$                            &     \\\hline
        \TS$\ww\LA n$ $\scr(n>5)$    & $\Rp+\ww{\su3}a+\ww{\su3}b+\Rm+\su{n-5}$    & $3$ \\
        \TS$(2^\sn3,1^\sn{n-5})$     & $\Lw+\End(\C^\sn3)$                         &     \\\hline
        \TS$\ww\LA n$ $\scr(n>3)$    & $\Rp+\Rm+\su{n-1}$                          & $4$ \\
        \TS$(3,1^\sn{n-2})$          & $\Lp2\Lm3+\Lw+(\Lm4+\Lm{-1})\Se\C^\sn{n-1}$ &     \\\hline
        \TS$\ww\LB 2$                & $\Rp+\Rm$                                   & $6$ \\
        \TS$(5)$                     & $(\Lw++\Lp{-1})\Lw-+\Lw+\Se\Ct$             &     \\\hline
        \TS$\ww\LB 3$                & $\Rp+\Rm+\su2$                              & $3$ \\
        \TS$(3,2^\sn2)$              & $\Lw++\Lm2+\Lw+\Lw-\Se\Ct$                  &     \\\hline
        \TS$\ww\LB n$ $\scr(n>3)$    & $\Rp+\Rm+\su2+\So n$                        & $4$ \\
        \TS $(3,2^\sn2,1^\sn{2n-6})$ & $\Lm2+\Lw+(\Lw-\Se\Ct+\Se\R^\sn n)$         &     \\\hline
        \TS$\ww\LB n$ $\scr(n>5)$    & $\R+\su6+\So{2n-11}$                        & $3$ \\
        \TS$(2^\sn6,1^\sn{2n-11})$   & $L\La^\sn2\C^\sn6$                          &     \\\hline        
        \TS$\ww\LC 3$                & $\R+\su3$                                   & $3$ \\
        \TS$(2^\sn3)$                & $L \St\C^\sn3$                              &     \\\hline        
        \TS$\ww\LC n$ $\scr(n>3)$    &$\R+\su3+\lie{sp}(n-3)$                      & $3$ \\
        \TS$(2^\sn3,1^\sn{2n-6})$    & $L \St\C^\sn3$                              &     \\\hline
        \TS$\ww\LD3$                 & $\Rp+\Rm+\su2$                              & $5$ \\
        \TS$(3^\sn2)$                & $\Lm8+(\Lw+ + \Lp{-1})\Lw-\Se\Ct$           &     \\\hline
        \TS$\ww\LD n$ $\scr(n>3)$    & $\Rp+\Rm+\su2+\So n$                        & $4$ \\
        \TS$(3,2^\sn2,1^\sn{2n-7})$  & $\Lm2+\Lw+(\Lw-\Se\Ct+\Se\R^\sn n)$         &     \\\hline        
        \TS$\ww\LD n$ $\scr(n>5)$    & $\R+\su6+\So{2n-12}$                        & $3$ \\
        \TS$(2^\sn{6},1^\sn{2n-12})$ & $L \La^\sn2\C^\sn6$                         &     \\\hline
        \hline\hline
      \end{tabular}
    }
  \caption[Beauville-bundle data of orbits of height three in classical orbits Lie algebras]{Beauville-data of orbits of
    height three in classical Lie algebras. $\lie n$ is decomposed as a $\LK$-module, $\End(\C^\sn3)$ is an
    $\ww{\SU(3)}a\ww{\SU(3)}b$-module under the action $(A,B)X\mapsto AXB^\sn{-1}$.}
  \label{Table:Big-table-classical}
\end{table}

\begin{table}[tb]
     {\small
      \begin{tabular}{lll}
        \hline\hline
        \TS Type$\backslash$Orbit & $\lie k$$\backslash$$\lie n$ & \makebox[0.5cm][r]{$\cohom\G\O$}\\
        \hline\hline
        \TS$\ww\LG2$             & $\R+\su2$                            & $6$ \\
        \TS$\scr 2\rthreea0$     & $L + L S^\sn3\Ct$                    &     \\\hline
        \TS$\ww\LF4$             & $\R+\su2+\su3$                       & $4$ \\
        \TS$\scr01\rtwoa00$      & $L^\sn3\Se\Ct+L^\sn2 \St\C^\sn3$     &     \\\hline
        \TS$\ww\LE6$             & $\R+\su2+\ww{\su3}a+\ww{\su3}b$      & $4$ \\
        \TS $\scr00\Dyatop0100$  & $L^\sn3\Se\Ct+L^\sn2\End(\C^\sn3)$   &     \\\hline
        \TS$\ww\LE7$             & $\R+\ww{\lie e}6$                    & $3$ \\
        \TS$\scr200\Dyatop0000$  & $L\ww V{27}$                         &     \\\hline
        \TS$\ww\LE7$             & $\R+\su2+\su6$                       & $4$ \\
        \TS$\scr000\Dyatop0010$  & $L^\sn3\Se\Ct+L^\sn2\La^\sn2\C^\sn6$ &     \\\hline
        \TS$\ww\LE8$             & $\R+\su2+\ww{\lie e}6$               & $4$ \\
        \TS$\scr0100\Dyatop0000$ & $L^\sn2\Se\Ct+L^\sn3\ww V{27}$       &     \\
        \hline\hline
      \end{tabular}
    }
  \caption[Beauville-bundle data of orbits of height three in exceptional Lie algebras]{Beauville-data of orbits of
    height three in exceptional Lie algebras. $\ww V{27}$ is a $27$-dimensional irreducible $\ww E6$-module.}
  \label{Table:Big-table-execptional}
\end{table}

\subsection{Examples of Standard Forms}
\label{sec:Examplex-of-Diagonalization}
With two exceptions, a nilpotent orbit of cohomogeneity three lies in a classical Lie algebra. It is either described by
the partition $(2^\sn3,1^\sn k)$ when the Lie algebra is of type $\LA$ or $\LC$, or by $(2^\sn6,1^\sn k)$ if the Lie
algebra has type $B$ or $D$. For these orbits, as the following examples show, it is fairly simple to `diagonalize'
elements of $\lie n$ via the action of $\LK$. This in fact only involves the spectral theorem applied to Hermitian
matrices.

 \begin{example}
  Consider the following action of the unitary group
  \begin{equation*}
    \begin{split}
      \LU(n)\times\Lie{Sym}(n,\C) & \to\Lie{Sym}(n,\C)\text{,}\\
      (g,Z)&\mapsto g\act Z:=g Z g^\sn t\text{.}
    \end{split}
  \end{equation*}
  F\mbox{}ix an element $Z$ in $\Lie{Sym}(n,\C)$. Then the Hermitian matrix $Z\Overline Z$ has real non-negative
  eigenvalues $\ww\mu1\geqs\dots\geqs\ww\mu n\geqs0$ and is diagonalized by some $g\in\LU(n)$ via conjugation,
  \begin{equation*}
    g(Z\Overline Z){\Overline g}^\sn t=\diag\left(\ww\mu1,\dots,\ww\mu n\right)=D\text{.}
  \end{equation*}
  Def\mbox{}ine the symmetric matrix $X=(\ww x{ij})$ by $X=g\act Z$ and observe that $X\Overline X=D$. Thus,
  \begin{gather*}
    \qq\sum{k=1}n |\ww x{ik}|^\sn2 = \ww\mu i\quad \text{and} \quad \qq\sum{k=1}n \ww x{ik} \ww{\Overline x}{kj}=0 \quad
    \text{for}\quad i\ne j\text{.}
  \end{gather*}
  Making use of the symmetry $\ww x{ij}=\ww x{ji}$ it follows that
  \begin{equation}
    \label{eq:sym-diag}
    \begin{split}
      \ww x{pq}\left(\ww\mu p-\ww\mu q\right)=&\ww\sum{j\ne q}\ww x{qp}\,|\ww x{pj}|^\sn2
      -\ww\sum{j\ne p}\ww x{pq}\,|\ww x{qj}|^\sn2+\left(\ww x{qp}\,|\ww x{pq}|^\sn2-\ww x{pq}\,|\ww x{qp}|^\sn2\right)\\
           =&\ww\sum{j\ne q}\left(\ww\sum{k}\ww x{qk}\,\ww{\Overline x}{kj}\right)\ww x{pj}
      - \ww\sum{j\ne p}\left(\ww\sum{k}\ww{x}{pk}\,\ww{\Overline x}{kj}\right)\ww x{qj}\\
      =&0\text{.}
     \end{split}
  \end{equation}
  If the eigenvalues are all dif\mbox{}ferent, $\ww\mu 1>\dots>\ww\mu n$, this implies that $X$ is diagonal. The general
  case follows immediately: let $M$ be the dense subset of $\Lie{Sym}(n,\C)$ whose elements have dif\mbox{}ferent
  eigenvalues and f\mbox{}ind a sequence $(\ww Xi)\subset\Lie M$ such that $\ww Xi\to X$. Then there exists elements $\ww
  gi\in\LU(n)$ such that $\ww gi\left(\ww Xi\ww{\Overline X}i\right)\qq{\Overline g}it$ is diagonal. By compactness
  we may assume that $\ww gi\to\ww g0\in\LU(n)$. Then $\ww g0\left(X\Overline X\right)\ww{\Overline g}0$ is diagonal
  with non-equal eigenvalues and~\eqref{eq:sym-diag} implies that $\ww g0\act X=(\ww g0 g)\act Z$ is diagonal.
\end{example}

\begin{example}
  In this example the symmetric matrices are replaced by skew-symmetric matrices,
  \begin{equation*}
    \begin{split}
      \LU(n)\times\La^\sn2\C^\sn n & \to\La^\sn2\C^\sn n\text{,}\\
      (g,Z)&\mapsto g\act Z:=g Z g^\sn t\text{.}
    \end{split}
  \end{equation*}
  Fix $Z\in\son n=\La^\sn2\C^\sn n=\List{Z\in\lie{gl}(n,\C)|Z+Z^\sn t=0}$ with $n$ even (the case for $n$ odd is
  similar).  The eigenvalues of the Hermitian matrix $Z\Overline Z$ are real and non-positive and there exists a
  $g\in\LU(n)$ with
  \begin{equation*}
    g(Z\Overline Z){\Overline g}^\sn t=\diag\left(\ww{\tilde\mu}1,\dots,\ww{\tilde\mu}n\right)=D \quad
    \text{where} \quad0\geqs\ww{\tilde\mu}1\geqs\dots\geqs\ww{\tilde\mu}n\text{.}
  \end{equation*}
  Put $X=g\act Z$ and assume $X$ to have $n$ dif\mbox{}ferent eigenvalues. Then $X$ is diagonalizable and $\C^\sn{n}$ is
  the direct sum of eigenspaces for both of the commuting matrices $X$ and $\Overline X$. As $X$ is skew-symmetric it
  follows that $\Overline Xv=-\Overline\la v$ whenever $v$ is a $\la$-eigenvector of $X$. Thus $v$ and $\Overline v$
  are eigenvectors of $X\Overline X$ corresponding to the same eigenvalue $-|\la|^\sn2$. Consequently,
  \begin{equation}
    \label{eq:skew-sym-diag-form} 
   X\Overline X=\diag\left(\ww\mu1,\ww\mu1,\ww\mu2,\ww\mu2,\dots,\ww\mu{n/2},\ww\mu{n/2}\right)
    ,\quad\ww\mu1>\dots>\ww\mu{n/2}\text{.}
  \end{equation}
  By a computation similar to~\eqref{eq:sym-diag} one gets that 
  \begin{equation*}
    \ww x{ij}(\ww\mu i-\ww\mu j)=0\text{,}
  \end{equation*}
  so $X=\ww g0\act Z$ has the block-diagonal form
  \begin{equation}
    \label{eq:skew-sym-anti-diag-form}
    \diag\left(\begin{smatrix}0&\ww x{12}\\-\ww x{12}&0\end{smatrix},\dots,\begin{smatrix}0&\ww x{n-1,n}\\-\ww x{n-1,n}&0
    \end{smatrix}\right)\text{.}
  \end{equation}
  Notice that the collection $\Lie M$ of elements in $\son n$ with $n$ dif\mbox{}ferent eigenvalues forms a dense subset
  of $\son n$. Thus, if X does not have $n$ dif\mbox{}ferent eigenvalues we may f\mbox{}ind a sequence $\ww Xi\in\Lie M$
  such that $\ww Xi\to X$. The previous argument implies that there exist $\ww hi\in\LU(n)$ such that $\ww hi\left(\ww
    Xi\ww{\Overline X}i\right) {\ww{\Overline h}i}^\sn t$ has the form in~\eqref{eq:skew-sym-diag-form} and $\ww hi\act
  X$ is given by (\ref{eq:skew-sym-anti-diag-form}) for all $i$. Since $\LU(n)$ is compact we may demand $\ww hi\to \ww
  h0\in\LU(n)$.  Now $(\ww h0 g)\act Z=\ww h0\act X$ is in the form given by~\eqref{eq:skew-sym-anti-diag-form}.
\end{example}  

\begin{example}
  Let $A\in\Lie M(n,\C)$. The self-adjoint matrix ${\Overline A}^\sn tA$ is diagonalized by some
  $q\in\SU(n)$, ${\Overline q}^\sn t({\Overline A}^\sn t A)q=\diag(\ww\la1,\dots,\ww\la n):=D$ and $\ww\la i\geqs0$. The
  $i$'th column $\ww qi$ of $q$ is an eigenvector of $A$ with eigenvalue $\ww \la i$. If we put $\ww
  ui=\qq\la i{-1/2}A\,\ww qi$ then $u:=\left(\ww u1\dots\ww un\right)\in\LU(n)$ and $A\,q=u\sqrt{D}$. Thus
  any $n\times n$ matrix can be diagonalized by the action
  \begin{equation*}
    \begin{split}
      \left(\LU(n)\times\SU(n)\right)\times\Lie M(n,\C)&\to\Lie M(n,\C)\text{,}\\
      \left((u,q),A\right)&\mapsto u\,A\,{\Overline q}^\sn t\text{.}
    \end{split}
  \end{equation*}
  Notice that this is actually the singular value decomposition.
\end{example}

\section{Potentials for HyperK\"ahler Metrics in Cohomogeneity Three}
In this chapter $\liegc$ is simple and $\O$ is a nilpotent orbit of cohomogeneity three under the action of the
compact group $\G$. Our aim is to f\mbox{}ind $\G$-invariant hyperK\"ahler metrics with K\"ahler potential on the complex
symplectic manifold $(\O,I,\Si)$ where $I$ is the natural complex structure and $\Si$ is the
Kirillov-Kostant-Souriau form.

Recall that a \emph{hyperK\"ahler potential} $\rho\colon \O\to\R$ is a function that is a simultaneous K\"ahler
potential for each of the complex structures, i.e. $\ww\om I = - i \ww\pa I \ww{\Overline\pa}I \rho$ etc. For our
purposes it is convenient to note that $i\ww\pa I\ww{\Overline\pa} I\rho=\tfrac12id(d-iId)\rho=\tfrac12 dId \rho$. In
general, the existence of a global K\"ahler potential is a very delicate matter. Indeed a compact manifold admits no
global K\"ahler potential. However, a local potential always exists (see eg.~\cite{GriffithsHarris}). There is no direct
analogue of this for hyperK\"ahler potentials. Generally, even the existence of a local hyperK\"ahler potential may not
be possible, see Swann~\cite{Swann:HK-quaternionicKahler}.

In this section we use the following abbreviation for the conjugation corresponding to the compact real form
\begin{equation*}
  Z'=\si Z,\quad Z\in\liegc\text{.}
\end{equation*}

\subsection{K\"ahler Potentials Depending on Three Invariants}
\label{sec:big-calculations}
Def\mbox{}ine three functions on the nilpotent orbit $\O$
\begin{gather*}
  \ww\eta1(X)=\inp{X}{X'},\quad  \ww\eta2(X)=-\inp{\DL{XX'}}{\DL{XX'}}\quad\text{and}\\
  \ww\eta3(X)=-\inp{\DL{XXX'}}{\DL{X'XX'}}, \quad X\in\O\text{,}
\end{gather*}
Observe that $\ww\eta2(X)=\ww\eta1(\DL{XX'})$ and $\ww\eta3(X)=\ww\eta1(\DL{XXX'})$, thus
$\ww\eta1,\ww\eta2,\ww\eta3:\O\to\R$ are positive $\G$-invariant functions.  Assume that $\rho:\O\to\R$ is a
$\G$-invariant K\"ahler potential for $(\O,I)$. The fact that the group $\G$ acts with cohomogeneity three allows us to
assume $\rho$ to be a function of our three invariant functions: $\rho=\rho( \ww\eta1,\ww\eta2, \ww\eta3)$. The K\"ahler
form is the second derivative of $\rho$, so by dif\mbox{}ferentiation one may f\mbox{}ind an explicit expression of
$\ww{\om}{I}$. To ease the notation let $\gZ\colon\liegc\to\liegc$ be the map
\begin{equation*}
  \gZ X=\DL{XX'X'X'X}+\DL{X'X'XXX'}+2\DL{X'XX'XX'}\text{.}
\end{equation*}
\begin{lemma}
  \label{lemma:omega-full}
  The K\"ahler form $\ww\om I$ def\mbox{}ined by $\rho$ is
    \begin{equation}
      \label{eq:omgaI-completly}
      \begin{split}
        \ww{\ww{&\om} I (\xia,\xib)}X=2\,\ww\rho1\im\inp\xia{\xibd}\\
        -4\,&\ww\rho2\im\inp\xia{2\DL{X'X\xibd}-\DL{XX'\xibd}}\\
        -2\,&\ww\rho3\im\bigl\{
        \!\!\!\begin{aligned}[t]
         &-2\,\inp\xia{\DL{XXX'X'\xibd}}+3\,\inp\xia{\DL{XX'XX'\xibd}}\\
         &+3\,\inp\xia{\DL{XX'X'X\xibd}}+3\,\inp\xia{\DL{X'XXX'\xibd}}\\
         &-12\,\inp\xia{\DL{X'XX'X\xibd}}+3\,\inp\xia{\DL{X'X'XX\xibd}}\bigr\}
        \end{aligned}
        \\
        +2\,&\ww\rho{11}\im\left(\inp{\xia}{X'}\inp{\xibd}X\right)\\
        -4\,&\ww\rho{12}
        \begin{aligned}[t]
          \im \bigl(
          &\inp{\xia}{X'} \inp{\xibd}{\DL{XX'X}}+\inp{\xia}{\DL{X'XX'}} \inp{\xibd}{X} \bigr)
        \end{aligned}
        \\
        +2\,&\ww\rho{13}\im\left(\inp\xia{X'}\inp{\xibd}{\gZ X'}+\inp\xia{\gZ X}\inp{\xibd}X\right)\\
        +8\,&\ww\rho{22}\im\left(\inp\xia{\DL{X'XX'}}\inp{\xibd}{\DL{XX'X}}\right)\\
        -4\,&\ww\rho{23}
        \begin{aligned}[t]
          \im\bigl(
          &\inp\xia{\DL{X'XX'}}\inp{\xibd}{\gZ X'}\\
          &+ \inp\xia{\gZ X}\inp{\xibd}{\DL{XX'X}}\bigr)
        \end{aligned}
        \\
        +2\,&\ww\rho{33}\im\left(\inp\xia{\gZ X}\inp{\xibd}{\gZ X'}\right)\text{,}
      \end{split}
    \end{equation}
  for any $\xia,\xib\in\Gamma(T\O)$ and $X\in\O$.
\end{lemma}

\begin{proof}
  We start by expressing $\ww\om I$ in terms of the derivatives of the invariant functions
  $\ww{\eta}{j}$,
  $j=1, 2, 3$,
  \begin{equation}
    \begin{split}
      \label{eq:omegaI-formelt}
      -2 \ww\om I
      &=dId\rho\\
      &=\ww\rho1\,dId\ww\eta1 + \ww\rho2\,dId\ww\eta2 + \ww\rho3\,dId\ww\eta3 \\
      &\quad+\ww\rho{11}\,d\ww\eta1\we Id\ww\eta1 +\ww\rho{22}\,d\ww\eta2\we Id\ww\eta2
      +\ww\rho{33}\,d\ww\eta3\we Id\ww\eta3\\
      &\quad+\ww\rho{12}\,(d\ww\eta1\we Id\ww\eta2 + d\ww\eta2\we Id\ww\eta1)\\
      &\quad+\ww\rho{13}\,(d\ww\eta1\we Id\ww\eta3 + d\ww\eta3\we Id\ww\eta1)\\
      &\quad+\ww\rho{23}\,(d\ww\eta2\we Id\ww\eta3 + d\ww\eta3\we Id\ww\eta2)\text{.}
    \end{split}
  \end{equation}
  The next step is to determine all dif\mbox{}ferentials $d\ww\eta j$ and $dId\ww\eta j$, $j=1,2,3$. However, only
  calculations involving $\ww\eta3$ will be carried out in detail, the other cases being similar if not simpler. Let us
  write $\ww\eta3=\ww\eta1\circ\phi$ where $\phi$ is the endomorphism of $\liegc$ def\mbox{}ined by $\phi(Z) =
  \DL{ZZZ'}$, $Z\in\liegc$. Dif\mbox{}ferentiating $\phi$ and $\ww\eta1$ is easy,
  \begin{eqnarray*}
    \ww{(d\ww\eta1)}{\phi(X)}(Z)&=&\inp{Z}{\DL{X'X' X}} + \inp{\DL{XXX'}}{Z'}\quad,\quad Z\in\lieg\\
    \ww{d\phi(\xia)}X&=&\DL{{\xia}XX'} + \DL{X{\xia}X'} + \DL{XX\xiad}\text{.}
  \end{eqnarray*}
  A repeated use of the Jacobi identity and the $\ad$-invariance of the Killing form give us the derivative of
  $\ww\eta3$
  \begin{equation*}
    \begin{split}
      \ww{d\ww\eta3(\xia)}{X}&=\inp\xia{\DL{\DL{XX'}X'X'X}-\DL{X'XX'X'X}+\DL{X'X'XXX'}}\\
      &\quad+\inp{\xiad}{\DL{XXX'X'X}+\DL{\DL{X'X}XXX'}-\DL{XX'XXX'}}\\
      &=2\re\qq\sum{j=1}3\inp\xia{\ww\gZ jX}\text{,}
    \end{split}
  \end{equation*}
  Here we have def\mbox{}ined $\ww\gZ1X=\DL{XX'X'X'X}$, $\ww\gZ2X=\DL{X'X'XXX'}$ and $\ww\gZ3X=2\,\DL{X'XX'XX'}$.
  Therefore $\ww{Id\ww\eta3(\xia)}{X}=2\im\ww\sum j\inp\xia{\ww\gZ jX}$, of which we are now able to f\mbox{}ind the
  derivative. However, this is somewhat more demanding. F\mbox{}irst observe that
  \begin{equation}
    \label{eq:dIdeta3-1}
    \begin{split}
      dId&\ww{\ww\eta3(\xia,\xib)}X\\
      &=\xia\ww|X (Y\mapsto Id\ww\eta3(\ww{(\xib)}{Y})) - \xib\ww|X (Y\mapsto Id\ww\eta3(\ww{(\xia)}{Y})) -
      \ww{Id\ww\eta3(\DL{{\xia}{\xib}})}X\\
      &=2\im\qq{\sum}{j=1}{3}
      \begin{aligned}[t]
        \bigl\{
        & \xia\ww|X ( Y\mapsto\inp{\DL{BY}}{\ww\gZ jY} ) - \xib\ww|X (Y\mapsto\inp{\DL{AY}}{\ww\gZ jY})\\
        &-\inp{\DL{XAB}}{\ww\gZ jX} \bigr\}.
      \end{aligned}
    \end{split}
  \end{equation}
  Letting $\ww\chi j$ denote the summands in~\eqref{eq:dIdeta3-1} involving $\ww\gZ j$ one f\mbox{}inds that
  \begin{equation*}
    \begin{split}
      \ww\chi1&=2\,\im\big\langle\xia\,,\,-\DL{{\xib}X'X'X'X}+\DL{\DL{{\xibd}X'}XXX'}-\DL{\DL{X{\xibd}X'}XX'}\\
      &+\DL{X'XXX'\xiad}+\DL{X'X'X'X\xib}-\DL{{\xib}X'X'X'X}-\DL{X{\xibd}X'X'X}\\
      &-\DL{XX'{\xibd}X'X}+\DL{XX'X'X\xibd}-\DL{XX'X'X'\xib}\,\big\rangle\text{.}
    \end{split}
  \end{equation*}
  For reasons we shall explain later, it is desirable to place the $\xib$'s in the above formula to the far right.
  Although the calculations are lengthy when using the Jacobi identity, it is a straight forward exercise. One
  f\mbox{}inds that
  \begin{equation*}
    \begin{split}
      \ww{\chi}1&=2\,\im\big\langle\xia\,,\,2\,\DL{XXX'X'\xibd}-6\,\DL{XX'XX'\xibd}+3\,\DL{XX'X'X\xibd}\\
      &-3\,\DL{XX'X'X'\xib}+3\,\DL{X'XXX'\xibd}+6\,\DL{X'XX'X'\xib}\\
      &-6\,\DL{X'X'XX'\xib}+3\,\DL{X'X'X'X\xib}\,\bigr\rangle\text{.}
    \end{split}
  \end{equation*}
  Similarly, the expressions for the remaining two summands are
  \begin{equation*}
    \begin{split}
      \ww{\chi}2&=2\im\bigl\langle\xia\,,\,-2\,\DL{XXX'X'\xibd}+4\,\DL{XX'XX'\xibd}-8\,\DL{X'XX'X\xibd}\\
      &+\DL{XX'X'X'\xib}+\DL{X'XXX'\xibd}+\DL{XX'X'X\xibd}-\DL{X'X'X'X\xib}\\
      &+2\,\DL{X'X'XX\xibd}+2\,\DL{X'X'XX'\xib}-2\,\DL{X'XX'X'\xib}\,\bigr\rangle
    \end{split}
  \end{equation*}
  and
  \begin{equation*}
    \begin{split}
      \ww{\chi}3&=-4\im\bigl\langle\xia\,,\,2\,\DL{XXX'X'\xibd}-4\,\DL{XX'XX'\xibd}+8\DL{X'XX'X\xibd}\\
      &-\DL{XX'X'X'\xib}-\DL{X'XXX'\xibd}-\DL{XX'X'X\xibd}-\DL{X'X'X'X\xib}\\
      &-2\,\DL{X'X'XX\xibd}-2\,\DL{X'X'XX'\xib}+2\,\DL{X'XX'X'\xib}\,\bigr\rangle\text{.}
    \end{split}
  \end{equation*}
  We are now in the position to write down the second derivative of $\ww\eta3$ in a somewhat shortened form,
  \begin{equation*}
    \begin{split}
      dId&\ww{\ww\eta3(\xia,\xib)}{X} = 12\im\bigl\langle\,\xia\,,\, -\tfrac23\,\DL{XXX'X'\xibd}-4\,\DL{X'XX'X\xibd}\\
        +&\DL{XX'X'X\xibd}+\DL{X'XXX'\xibd}+\DL{XX'XX'\xibd}+\DL{X'X'XX\xibd}\,\bigr\rangle\text{.}
    \end{split}
  \end{equation*}
  As already noted, the determination of $d\ww\eta j$, $dId\ww\eta j$, $j=1,2$ represents a similar exercise, moreover
  it is already present i ~\cite{KobakSwann:HK-pot-in-cohom-two}
  \begin{gather*}
    \ww{d\ww\eta1(\xia)}X = 2\re\inp\xia{X'},\quad \ww{dId\ww\eta1(\xia,\xib)}X = -4\im\inp\xia\xibd,\\
    \ww{d\ww\eta2(\xia)}X = -4\re\inp\xia{\DL{X'XX'}} \quad \text{and}\\
    \ww{dId\ww\eta2(\xia,\xib)}X = 8\im\inp\xia{2\,\DL{X'X\xibd}-\DL{XX'\xibd}}.
  \end{gather*}
  To f\mbox{}inish, note that with $\gZ X=\qq{\sum}{j=1}{3}\ww\gZ jX$ we have
  \begin{equation*}
    \begin{split}
      \ww{(d\ww\eta1\we& Id\ww\eta3 + d\ww\eta3\we Id\ww\eta1)(\xia,\xib)}{X}\\
      &=
      \begin{aligned}[t]
        &4\bigl(\re\inp\xia{X'}\im\inp\xib{\gZ X}-\re\inp\xib{X'}\im\inp\xia{\gZ X}\\
        &+\re\inp\xia{\gZ X}\im\inp\xib{X'}-\re\inp\xib{\gZ X}\im\inp\xia{X'} \bigr)
      \end{aligned}
      \\
      &=-4\im\left(\inp\xia{X'}\inp{\xibd}{\gZ X'}+\inp\xia{\gZ X}\inp{\xibd}X\right),
    \end{split}
  \end{equation*}
  with similar expressions for the other summands of~\eqref{eq:omegaI-formelt}. Collecting the
  relevant formulæ completes the proof.
\end{proof}
Supposing $g(\cdot,\cdot)=\ww\om I(I\cdot,\cdot)$ to be non-degenerate we may def\mbox{}ine an endomorphism $J$ of the
tangent bundle $T\O$ by
\begin{equation}
  \label{eq:J-def}
  g(\xia,\xib)=\re \Si(J\xia,\xib)
\end{equation}

\begin{lemma}
  The endomorphism $J$ is given by
  \begin{equation}
    \label{eq:J-cohom3}
    \begin{split}
      \ww{(J& \xia)}X =-2\ww\rho1 \DL{X\xiad}\\
      +&4\, \ww\rho2 \left( 2\,\DL{XX'X\xiad}-\DL{XXX'\xiad} \right)\\
      +&2\, \ww\rho3\bigl(
      \begin{aligned}[t]
        & -2\,\DL{XXXX'X'\xiad}+3\,\DL{XXX'XX'\xiad}\\
        & +3\,\DL{XXX'X'X\xiad}+3\,\DL{XX'X'XX\xiad}\\
        &-12\,\DL{XX'XX'X\xiad}+3\,\DL{XX'XXX'\xiad}\,\bigr)
      \end{aligned}
      \\
      -&2\, \ww\rho{11} \inp\xiad X \DL{XX'}\\
      +&4\, \ww\rho{12}
      \begin{aligned}[t]
        \bigl(
        &  \inp\xiad X \DL{XX'XX'}+ \inp\xiad{\DL{XX'X}} \DL{XX'} \bigr)
      \end{aligned}
      \\
      -&2\, \ww\rho{13}\bigl(\inp\xiad X \DL{X(\gZ X)}+\inp{\xiad}{\gZ X'}\DL{XX'}\bigr)\\
      -&8\, \ww\rho{22} \inp\xiad{\DL{XX'X}} \DL{XX'XX'}\\
      +&4\, \ww\rho{23}
      \begin{aligned}[t]
        \bigl(
        &  \inp\xiad{\DL{XX'X}} \DL{X(\gZ X)} + \inp\xiad{\gZ X'} \DL{XX'XX'} \bigr)
      \end{aligned}
      \\
      -&2\, \ww\rho{33} \inp\xiad{\gZ X'}\DL{X(\gZ X)}\text{,}
    \end{split}
  \end{equation}
  for any $\xia,\xib\in\Gamma(T\O)$ and $X\in\O$.
\end{lemma}
\begin{proof}
  $J$ is determined by the relation
  \begin{equation*}
    -\re\inp{\ww J X\xia}{B}=\ww{\ww\om I(I \xia,\xib)}{X}\text{.}
  \end{equation*}
  and the non-degeneracy of the form $\inp\cdot\cdot$. As $I\xia=i\,\xia$, the right-hand side above equals the
  right-hand side of equation~\eqref{eq:omgaI-completly} with the slight modif\mbox{}ication of replacing `$\im$' with
  `$\re$' throughout. The result can now be read of\mbox{}f equation~\eqref{eq:omgaI-completly}.
\end{proof}

\begin{remark}
  \label{Remark:Jxia}
  Notice that $J\xia$ remains in the subalgebra generated by $\List{X,X',A,A'}$.
\end{remark}

\subsection{Generic Cohomogeneity-Three Orbits}
\label{sec:general_cohom_three_orbits}
For any $X\in\O$, assume the Lie span of $\List{X,X'}$ to be embedded into a $\si$-invariant subalgebra isomorphic to
three copies of $\slc$; we shall denote such an algebra by $\slT$.

\begin{remark}
  \label{Remark:general-cohom-three}
  The majority of cohomogeneity-three orbits satisfy the above assumption, and as such will be referred to as
  \emph{generic cohomogeneity three orbits}. In fact it shall become clear that $\son7$ posses the only \emph{special
    cohomogeneity three orbit} (see \S~\ref{section:Rel-w-Hermitian-Sym-Spaces}). The orbits of cohomogeneity three are
  listed in Table~\ref{Table:All-cohom-three}.
\end{remark}

\begin{table}[h]
  \centering
  \begin{tabular}[t]{cc}
    \multicolumn{2}{c}{\textbf{Generic}}\\
    \hline\hline
    \TS Type & Orbit\\
    \hline
    \TS $\ww\LA n$ & $(2^\sn3,1^\sn{n-5})$\\
    \TS $\ww\LB{(n-1)/2}$, $\ww\LD{n/2}$ & $(2^\sn6,1^\sn{n-12})$\\
    \TS $\ww\LC n$ & $(2^\sn3,1^\sn{2n-6})$\\
    \TS $\ww\LE7$ & $\scr 200\Dyatop0000$\\
    \hline\hline
  \end{tabular}
  \qquad
  \begin{tabular}[t]{cc}
    \multicolumn{2}{c}{\textbf{Special}}\\
    \hline\hline
    \TS Type & Orbit\\
    \hline
    \TS $\ww\LB3$ & $(3,2^\sn2)$\\
    \hline\hline
  \end{tabular}
  \caption[Orbits of cohomogeneity three in simple Lie algebras]{Orbits of cohomogeneity three in simple Lie algebras}
  \label{Table:All-cohom-three}
\end{table}

The $\si$-invariance of $\slT$ tell us that the restriction $\tilde\si:=\si\ww{|}{\slc^\sn3}$ is a compact real
form of $\slT$, so the $+1$ eigenspace of $\tilde\si$ is isomorphic to three copies of $\su2$,
\begin{equation*}
  (\slc^\sn3)^\sn{\tilde\si}=\ww{\su2}+\op\ww{\su2}0\op\ww{\su2}-\text{.}
\end{equation*}
Let $\ww\slc\de$ denote the complexf\mbox{}ication of $\ww{\su2}{\de}$, where $\de$ refers to any of the three
subscripts $\List{+,0,-}$. Thus
\begin{equation*}
  X\in\slc^\sn3=\ww\slc+\op\ww\slc0\op\ww\slc-\subset\liegc\text{.}
\end{equation*}
We write $X=\ww X+ + \ww X0 + \ww X-$, $\ww X\de\in\ww\slc\de$ and assume $X$ to be generic (i.e. each $\ww X\de$ is
nonzero). Notice that $\O\cap\ww\slc\de$ is the minimal orbit of $\ww\slc\de$ and is of cohomogeneity one under the
action of $\ww{\SU(2)}{\de}$. It follows from Remark~\ref{remark:R-times-G} that
\begin{equation}
  \label{eq:aXbXcX}
  \ww Z{abc}=a\,\ww X+ + b\,\ww X0 + c\,\ww X-\in\O
\end{equation}
for all $a,b,c\in\Rp$. The element $Z=a\,\ww X+ + b\,\ww X0$ is the limit of~\eqref{eq:aXbXcX} as $c$ approaches zero so
it lies in the closure of $\O$. Considered as endomorphisms of $\slT$ the rank of $\ww\ad Z$ is obviously strictly less
than the rank of $\ad_{Z_{abc}}$, since the three subalgebras commute. On the other hand, the orthogonal complement of
$\slT$ is invariant under the adjoint action of $\slT$. Since $\ww Z{abc}$ is an element of $\O$ the rank the of
$\ad_{Z_{abc}}$ equals the constant rank of $\ww\ad X$ on $\slT$ and its orthogonal complement for all $a,b,c\in\Rp$.
By continuity the rank of $\ww\ad Z$ cannot exceed the rank of $\ad_{Z_{abc}}$ on the orthogonal
complement. Taken together, $Z$ cannot be an element of $\O$, and via~\eqref{eq:Boundary-relation} therefore lies in an
orbit of height $2$ or $1$.  Repeating this argument we f\mbox{}ind that $\ww X\de$ lies in the minimal orbit of $\liegc$
and is conjugate to a highest root vector.

Expressing everything in terms of $\G$-invariant objects, we may restrict our attention
to points of the form
\begin{equation}
  \label{eq:decomp-X}
  X=\ww X+ +\ww X0 + \ww X- = s\,\ww E+ + r\,\ww E0 + t\,\ww E-\text{,}
\end{equation}
where $s,r,t>0$ and each $\ww E\de$ is a highest root vector in $\liegc$. Choose $\ww F\de=-\siEd\de$, $\ww H\de=\DL{\ww
  E\de\ww F\de}$, so that $\List{\ww H\de,\ww E\de,\ww F\de}$ is a standard basis of $\ww\slc\de$. Let
$\inp[\sn\de]\act\act$ denote the negative of the Killing form of $\ww\slc\de$. By Schur's Lemma we have
$\inp\act\act\ww|{\ww\slc\de} = \qq k\de2\inp[\sn\de]\act\act$, where the constant $\qq k\de2$ is strictly positive
because $\ww{\su2}\de\subset\lieg$. The algebras $\ww\slc\de$ correspond each to a highest root in the same semisimple
Lie algebra, so they are conjugate by the Weyl group and we have $\qq k\de2=k^\sn2$.

Our three $\G$-invariant functions on $\O$ are given by
\begin{equation*}
  \ww\eta i=2^\sn{i+1} k^\sn2\left(s^\sn{2i}+r^\sn{2i}+t^\sn{2i}\right),\quad i=1,2,3\text{.}
\end{equation*}

Let us use the abbreviations $\ww\rho i = \pa\rho /\pa\ww\eta i$ and $\ww\rho u = \pa\rho/\pa u$, where $i$ is one of
$1,2,3$ and $u$ ranges in $\List{r,s,t}$. Then the derivative $d \rho =\ww\rho1 d\ww\eta1 +\ww\rho2 d\ww\eta2 + \ww\rho3
d\ww\eta3$ gives rise to the following system of equations
\begin{equation}
  \label{eq:rho_coeff}
  \begin{pmatrix}
    \ww\rho s\\
    \ww\rho r\\
    \ww\rho t
  \end{pmatrix}
  =8 k^\sn2
  \begin{pmatrix}
    s&4s^\sn3&12s^\sn5\\
    r&4r^\sn3&12r^\sn5\\
    t&4t^\sn3&12t^\sn5\\
  \end{pmatrix}
  \begin{pmatrix}
    \ww\rho 1\\
    \ww\rho 2\\
    \ww\rho 3
  \end{pmatrix}\text{.}
\end{equation}

\begin{remark}
  The constant $k^\sn2$ only depends on the Lie algebra $\liegc$. If we consider the highest root decomposition
  \begin{equation*}
    \liegc=\ww\slc+ \op \liekc \op \qq S+1\ot\ww V+\text{,}
  \end{equation*}
  we notice that $\inp{\ww E+ }{\siEd+} = 4+\ww\dim\C\ww
  V+$. On the other hand $\inp{\ww E+}{\siEd+}=k^\sn2\inp[\sn+]{\ww E+}{\siEd+}=4k^\sn2$, so one may use
  equation~\eqref{eq:dim-of-V} to f\mbox{}ind the value of $k^\sn2$.
  \begin{table}[h]
    \centering
    \newcommand{\TSs}{\vrule height 14pt depth 6pt width 0pt}
    \leavevmode
    \begin{tabular}[t]{c|ccccccc}
      \hline\hline
      \TSs Type &$\ww\LA n,\ww\LC n$ & $\ww\LB{(n-1)/2},\ww\LD{n/2}$ & $\ww\LG2$ &
      $\ww\LF4$ & $\ww\LE6$ & $\ww\LE7$ & $\ww\LE8$\\
      \hline
      \TSs $k^\sn2$ & $\tfrac{n+1}2$ & $\tfrac{n-2}2$ & $2$ & $\tfrac92$ & $6$ & $9$ & $15$\\
      \hline\hline
    \end{tabular}
    \caption{The constant $k^\sn2$.}
  \end{table}
\end{remark}

\subsubsection{The Regular Orbit of $\slT$}
\label{sec:restrict-to-threesl2C}
The possible choices for $\rho$ are narrowed by imposing the endomorphism $J$ in~\eqref{eq:J-def} to be an almost
complex structure. For the moment we restrict our attention to the subalgebra $\slT$, and
in fact we are considering the regular orbit $\ww\O{reg}=\ww\O+\times\ww\O0\times\ww\O-$ of $\slT$, where $\ww\O\de$ is
the nonzero nilpotent orbit in $\ww\slc\de$. By Remark~\ref{Remark:Jxia}, $J$ is indeed an endomorphism of $\slT$. Let
$A\in\slc^\sn3$ and put $\qq\xi A\de=\DL{A\ww X\de}$. In Lemma~\ref{lemma:omega-full} we expressed the K\"ahler form
$\ww\om I$ in terms of the functions $\ww\eta i$, $i=1,2,3$. We now use the functions $\ww\eta\de$ def\mbox{}ined by
\begin{equation*}
  \ww\eta\de =\ww\eta1(\ww X\de)=k^\sn2\inpsl{\ww X\de}{\siXd\de}\text{,}
\end{equation*}
where $\Sl$ is short for $\slc$. Notice that $\ww\eta+=4\,k^\sn2s^\sn2$ etc.

The K\"ahler form can be written as
\begin{equation*}
  \begin{split}
    \ww{\ww\om I(\xia,&\xib)}X =2k^\sn2\im\bigl\{\\
    &
    \begin{aligned}[t]
      &\ww\rho+\inpsl{\qq\xi A+}{\si\qq\xi B+} +
      \ww\rho0\inpsl{\qq\xi A0}{\si\qq\xi B0} + \ww\rho-\inpsl{\qq\xi A-}{\si\qq\xi B-}\\
      &+\ww\rho{++}k^\sn2\inpsl{\qq\xi A0}{\siXd+}\inpsl{\si\qq\xi B+}{\ww X+}\\
      &+\ww\rho{+0}k^\sn2 \left(\inpsl{\qq\xi A+}{\siXd+}\inpsl{\si\qq \xi B0}{\ww X0} +
        \inpsl{\qq\xi A0}{\siXd0}\inpsl{\si\qq\xi B+}{\ww X+} \right) \\
      &+\ww\rho{+-}k^\sn2 \left( \inpsl{\qq\xi A+}{\siXd+}\inpsl{\si\qq\xi B-}{\ww X-} +
        \inpsl{\qq\xi A-}{\siXd-}\inpsl{\si\qq\xi B+}{\ww X+} \right)\\
      &+\ww\rho{00}k^\sn2\inpsl{\qq\xi A0}{\siXd0}\inpsl{\si\qq\xi B0}{\ww X0}\\
      &+\ww\rho{0-}k^\sn2 \left( \inpsl{\qq\xi A0}{\siXd0}\inpsl{\si\qq\xi B-}{\ww X-} +
        \inpsl{\qq\xi A-}{\siXd-}\inpsl{\si\qq\xi B0}{\ww X0} \right)\\
      &+\ww\rho{--}k^\sn2\inpsl{\qq\xi A-}{\si \ww X-}\inpsl{\si\qq\xi B-}{\ww X-} \bigr\}\text{,}
  \end{aligned}  
\end{split}
\end{equation*}
where $\ww\rho{++} = \pa\ww\rho+ /\pa\ww\eta+$, etc. The endomorphism $J$ is given by
\begin{equation*}
  \begin{split}
    \ww JX(\qq\xi A+)=
    -2\,\ww\rho+&\DL{\ww X+\si\qq\xi A+}-2\,k^\sn2\inpsl{\si\qq\xi A+}{\ww X+}\bigl(\ww\rho{++}\,\DL{\ww X+\siXd+}\\
    &+ \ww\rho{+0}\,\DL{\ww X0\siXd0} + \ww\rho{+-}\,\DL{\ww X-\siXd-}\bigr)
  \end{split}
\end{equation*}
and similarly for $\qq\xi A0$, $\qq\xi A-$. Evaluating at $\ww H+$ and $\ww E+$ we f\mbox{}ind that $\ww JX\ww
H+=-4\,s\,\ww\rho+\,\ww E+$ and
\begin{gather*}
  \ww JX\ww E+=2\,s\left(\ww\rho ++\ww\eta+\,\ww\rho{++}\right)\ww H+ +
  2\,\ww\eta+\,s^\sn{-1}\left(r^\sn2\,\ww\rho{+0}\,\ww H0 + t^\sn2\ww\rho{+-}\,\ww H- \right)\text{.}
\end{gather*}
If $J^\sn2$ is to be $-1$, then the $\ww\slc0\op\ww\slc-$-component of the following must vanish,
\begin{equation}
  \label{eq:Jsqrt-on-H}
    \qq JX2\ww H+=-8\,s^\sn2\,\ww\rho+\left(\ww\rho+ + \ww\eta + \,\ww\rho{++}\right)\ww H+ - 
    8\,\ww\rho+\,\ww\eta+\left(r^\sn2\ww\rho{+0}\,\ww H0 + t^\sn2\,\ww\rho{+-}\,\ww H-\right)\text{.}
\end{equation}
Thus $0=2\,\ww\rho+\,\ww\rho{+0}=\pa(\qq\rho+2)/\pa\ww\eta0$, so $\ww\rho{+0}=0$, and similarly $\ww\rho{+-}=0$ We
conclude that $J^\sn2=-1$ on $\slT$ if and only if $\ww\rho{+0}=\ww\rho{+-}=\ww\rho{0-}=0$, in which case each of
the $\slc$-summands is preserved.

Consequently, the condition $J^\sn2=-1$ is equivalent to
\begin{equation*}
  k^\sn2=2\,\ww\eta\de\,\ww\rho\de\left(\ww\rho\de + \ww\eta\de\,\ww{\rho}{\de\de}\right) =
  \frac\pa{\pa\ww\eta\de}\left(\ww\eta\de\,\ww\rho\de\right)^\sn2 
\end{equation*}
for each value of $\de$. Hence,
\begin{equation}
\label{eq:rho-delta}
  \qq\rho\de2=\left(k^\sn2\,\ww\eta\de+\frac{\ww c\de}4\right)/\qq \eta\de2\text{,}
\end{equation}
for some constant $\ww c\de$. We demand that $\ww c\de\geqslant0$ in order for $\ww\rho\de$ to be def\mbox{}ined for all
$\ww X\de$. Notice that in terms of the parameter $s$, we get
\begin{equation}
  \label{eq:rho-s-by-s}
  \qq\rho s2=16\,k^\sn4 + \ww c+/s^\sn2\text{,}
\end{equation}
and similarly for $\ww\rho r$ and $\ww\rho t$.

\begin{remark}
  \label{Remark:unique_hKpot}
  If we had considered the case $\liegc=\slc$ with inner product $k^\sn2\inp\act\act$ and used the methods described in
  Section~\ref{sec:big-calculations} with one invariant function, $\eta=\ww\eta1$, then the potential of the nonzero
  nilpotent orbit in $\slc$ would have been given by equation~\eqref{eq:rho-delta} and the corresponding metric would
  have looked like
\begin{equation}
    \label{eq:explicit-cohom-one-metric}
  \newcommand{\TSs}{\vrule height 10pt depth 0pt width 0pt}
  \leavevmode
  \begin{split}
    \frac{\ww{g(\xia,\xib)}X}{{\sqrt{4\,k^\sn2\,\eta+c}}} = \frac{k^\sn2}{\eta^\sn2}\re
    \left(\eta\inp\xia\xibd-\frac{2\,k^\sn2\eta\,+c}{4\,k^\sn2\,\eta+c}\inp\xia{X'}\inp{\xibd}X\right)\text{.}
  \end{split}
\end{equation}
Moreover, the potential is a hyperK\"ahler potential if and only if $\ww c\de=0$
(see~\cite{KobakSwann:HK-pot-in-cohom-two}).
\end{remark}

\begin{proposition}
  An $\SU(2)^\sn3$-invariant hyperK\"ahler structure of the principal orbit of $\slT$ that admits a K\"ahler
  potential and has the Kirillov-Kostant-Souriau complex symplectic form, is a product of $\SU(2)$-invariant structures
  (given in~\eqref{eq:explicit-cohom-one-metric}) on each factor. \qed
\end{proposition}

\subsubsection{Other Lie Algebras}
Recall that we are considering a nilpotent orbit $\O\subset\liegc$ such that $X\in\O$ is decomposed via
equation~\eqref{eq:decomp-X}. In~\S\ref{sec:restrict-to-threesl2C} we described the action of $J$ on tangent vectors
$\xia\in\slc^\sn3$ and found that the potential must satisfy equation~\eqref{eq:rho-s-by-s} for $J$ to be almost
complex. To get a complete picture we need to consider the case when $\xia$ lies in the Killing-orthogonal complement of
$\slT$. Additional obstructions in choosing $\rho$ may emerge this way. As discussed in
Section~\ref{Section:StandardTriple} $\ww X\de$ is a nilpotent element of $\liegc$ and via
\eqref{eq:Boundary-relation} lie in the minimal orbit of $\liegc$. By Proposition~\ref{prop:DancerSwann-cohom1} we have
a Killing-orthogonal decomposition corresponding to a highest root,
\begin{equation}
  \label{eq:highest-root-decom}
  \liegc\cong\ww\slc\de\op\qq{\lie k}\de\C\op(\qq S\de1\ot\ww V\de)\text{.}
\end{equation}
Recall that $\ww V\de$ is a non-trivial $\qq{\lie k}\de\C$-module and that $\qq{\lie k}\de\C$ is the centralizer of
$\ww\slc\de$. Note that $\ww\slc-\subset\qq{\lie k}+\C$, so, as an $\ww\slc-$-module, $\ww V+$ may consist of trivial
and fundamental $\ww\slc-$-modules; similarly for $\ww\slc0$. As we will see below, it suf\mbox{}f\mbox{}ices to
consider the following two situations:

\medskip

\noindent\textup{(i)} Assume the existence of a trivial $\ww\slc0\op\ww\slc-$-module $\C^\sn r$ in $\ww V+$. Choose the
dimension $r$ to be maximal, so that the real structure $\si$ preserves the module $\qq S+1\ot\C^\sn r$. We now describe
the action of $\si$.  Using the identif\mbox{}ication $\qq S+1\ot\C^\sn r=\ww S1\op\dots\op\ww Sr$, $\ww Sj \cong \qq
S+1$, let $\ww\si{ik}$ be the component of $\si$ that sends $\ww Si$ to $\ww Sk$. Note that $j \circ \si$
($j\in\H\cong\qq S+1$) commutes with the $\slc$-action and by Schur's Lemma is a scalar, hence $\ww\si{i k}=\ww\mu{i k}
j$, $\ww{\mu}{i k}\in\C$.  Since $\si$ is an involution, $(\ww\mu{ik})(\Overline{\ww\mu{ik}})^t=-1$ and $\C^\sn r$
admits a quaternionic structure, $\lie j\colon\ww ei \mapsto \qq\sum{k=1}r\ww\mu{ik}\ww ek$.  Thus $\si$ acts on $\qq
S+1\ot\C^\sn r$ as $j\ot \lie j$.
  
To determine $J$ on the module $\qq S+1\ot\C^\sn r$ we choose a basis on which $\ww E+$ acts as
$\begin{smatrix}0&1\\0&0\end{smatrix}$. In this picture any tangent vector $\xia=\DL{AX}\in \qq S+1\ot\C^\sn r$ has the
form ${1 \choose 0}\ot w$. We now see the convenience of the positions of `$\xia$' in~\eqref{eq:J-cohom3}.
We immediately get
\begin{equation*}
  \ww JX\xia=-2\,s \left( \ww\rho1+4\,s^\sn2\ww\rho2+12\,s^\sn4\ww\rho3 \right)\begin{smatrix}0\\1\end{smatrix}\ot\lie jw
    =-4\,k^\sn2\ww\rho s\begin{smatrix}0\\1\end{smatrix} \ot\lie jw\text{.}
\end{equation*}
Hence, as an endomorphism of $\qq S+1 \ot W$, $J^\sn2$ is $-1$ if and only if $\qq\rho s2=16 k^\sn4$. Comparing with
equation~\eqref{eq:rho-s-by-s} we now know that $\ww c+=0$ if $\ww V+$ has a trivial $\ww\slc0\op\ww\slc-$-module.

\medskip

\noindent\textup{(ii)} Next, consider the case of a tangent vector $\xia$ lying in a $\slT$-module $\qq S+1\ot\qq
S-1$. This is a submodule of $\liegc$ contained in the Killing orthogonal complement of $\slT$. It is convenient to
choose a basis such that $\ww E\pm$ acts on $\qq S\pm1$ as $\begin{smatrix}0&1\\0&0\end{smatrix}$ and $\si=j \ot j$.
In other words, $X$ acts as
\begin{equation}
  \label{eq:adX}
  s \begin{smatrix}0&1\\0&0\end{smatrix} \ot\Id + t\,\Id\ot\begin{smatrix}0&1\\0&0\end{smatrix}\text{.}
\end{equation}
We choose two independent vectors that span the image of $\ww\ad X$
\begin{equation*}
  \ww\xi1 = \begin{smatrix}1\\0\end{smatrix} \ot \begin{smatrix}1\\0\end{smatrix} \quad\text{and}\quad
  \ww\xi2 = s\begin{smatrix}1\\0\end{smatrix} \ot \begin{smatrix}0\\1\end{smatrix} + t \begin{smatrix}0\\1\end{smatrix}
  \ot \begin{smatrix}1\\0\end{smatrix}\text{.}
\end{equation*}
Now if $A\in\qq S+1\ot\qq S-1$ then $\xia\ww|X=\DL{AX}$ is in the span of $\ww\xi1$ and $\ww\xi2$ and
using~\eqref{eq:adX} the action of $X$ and $X'$ on these vectors is easily found to be
\begin{gather*}
  \DL{X\ww\xi1} = 0,\quad\DL{X'\ww\xi1} =\xitd,\quad\DL{X\ww\xi2} = 2\,s\,t\,\ww\xi1,\quad\DL{X'\ww\xi2} =
  -\left(s^\sn2+t^\sn2\right)\xiod\text{.}
\end{gather*}
Using~\eqref{eq:J-cohom3} we immediately get
\begin{gather*}
  \ww JX\ww\xi1=-2\left(\ww\rho1 + 4\,\ww\rho2\left(s^\sn2+t^\sn2\right) +
    12\,\ww\rho3\left(s^\sn4+s^\sn2t^\sn2+t^\sn4\right)\right)\ww\xi2\\
  \ww JX\ww\xi2 = 2\left(\ww\rho1\left(s^\sn2+t^\sn2\right) + 4\,\ww\rho2\left(s^\sn4+t^\sn4\right) +
    12\,\ww\rho3\left(s^\sn6+t^\sn6\right)\right)\ww\xi1\text{,}
\end{gather*}
which via~\eqref{eq:rho_coeff} may be expressed as
\begin{gather*}
  \ww JX\ww\xi1=-\frac{s \ww\rho s-t\ww\rho t}{4k^\sn2\left(s^\sn2-t^\sn2\right)}\ww\xi2\quad \text{and} \quad \ww
  JX\ww\xi2 = \frac{s\ww\rho s+t\ww\rho t}{4 k^\sn2}\ww\xi1\text{.}
\end{gather*}
Hence,
\begin{equation*}
  J^\sn2\ww\xi1=-\frac{s^\sn2\qq\rho s2-t^\sn2\qq\rho t2}{16\,k^\sn4\left(s^\sn2-t^\sn2\right)}\ww\xi1
\end{equation*}
and via~\eqref{eq:rho-s-by-s} we see that $J^\sn2$ is $-1$ on the module $\qq S+1\ot\qq S-1$ if and only if $\ww c+=\ww
c-$.

\medskip

With these two considerations at hand we are ready to draw the conclusions. By Table~\ref{Table:All-cohom-three} there
are four orbits to consider: $(2^\sn3,1^\sn{n-6})$ in $\sln n$, $(2^\sn6,1^\sn{n-12})$ in $\son n$,
$(2^\sn3,1^\sn{2n-6})$ in $\spn n$ and $\scr 200\Dyatop0000$ in $\qq{\lie e}7\C$.  We thus need to decompose $\ww V+$
under the action of $\ww\slc0\op\ww\slc-$ when $\G$ is a classical Lie group or the exceptional Lie Group $\ww\LE7$.
These representations can be found in~\cite{KobakSwann:HK-pot-in-cohom-two}, and are is listed in
Table~\ref{Table:Wolf-data} for convenience. Notice that in the case $\G=\SO(n)$ the centralizer is $\qq{\lie
  k}+\C=\slc\op\son{n-4}$, and a priori there are two choices of $\ww\slc0\op\ww\slc-\subset\liekc$. One consists in
taking $\ww\slc0=\slc$ and $\ww\slc-\subset\son{n-4}$ and the other $\ww\slc0\op\ww\slc-\subset\son{n-4}$.  However, the
f\mbox{}irst one is not an option because it would allow the appearance of copy of $\qq S+1\ot\qq S01\ot\qq S-1$ and
Jordan normal form would therefore be incorrect.
\begin{table}[htb]
  \begin{center}
    \begin{tabular}{ccc}
      \hline\hline
      \TS$\G$                  & $\qq{\lie k}+\C$   & $\ww V+$\\
      \hline
      \TS$\SP(n)\; (n\geqs3)$  & $\spn{n-1}$        & $\qq S01\op\qq S-1\op\C^\sn{2n-6}$\\
      \TS$\SU(n)\; (n\geqs6)$  & $\C\op\sln{n-2}$   & $\qq S01\op\qq S-1\op\C^\sn{n-6}$ \\
      \TS$\SO(n)\; (n\geqs12)$ & $\slc\op\son{n-4}$ & $\Ct\ot(\qq S01\op\qq S-1\op\R^\sn{n-12})$\\
      \TS$\ww\LE7$             & $\son{12}$         & $\C^\sn8\ot(\qq S01\op\qq S-1)$\\
      \hline\hline
    \end{tabular}
    \caption[The centralizer $\lie k_{\sscr +}^{\sscr \C}$ for the classical groups and $\Lie E_{\sscr 7}$]{The
      centralizer $\lie k_{\sscr +}^{\sscr \C}$ for classical groups and $\Lie E_{\sscr 7}$. The module $V_{\sscr +}$ is
      decomposed under the action of $\slc_{\sscr 0}\op\slc_{\sscr -}$.}
    \label{Table:Wolf-data}
  \end{center}
\end{table}
We conclude that $\liegc$ always contains copies of $\qq S+1\ot\qq S01$ and $\qq S+1\ot\qq S-1$, that copies of $\qq
S+1\ot\qq S01\ot\qq S-1$ never occur, and that $\ww V+$ contains a trivial $\ww\slc0\op\ww\slc-$-module unless $\lieg$
is one of $\lie{sp}(3)$, $\su6$, $\So{12}$ or $\ww{\lie e}7$. In particular the constants $\ww c\de$ of the
potential coincide,
\begin{equation*}
  c:=\ww c+=\ww c0=\ww c-\text{,}
\end{equation*}
and if $\ww V+$ does not have a trivial summand we get a one-parameter family of hyperK\"ahler metrics
with K\"ahler potential. Alternatively the constant $c$ must be zero and we get a unique potential.

\subsubsection{An Explicit Potential}
\label{subsec:explicit-pot}
If $\ww c+ = \ww c0 = \ww c-=0$ the potential is $\rho=4\,k^\sn2\left(s+t+r\right)$ up to addition of some constant.
This can be expressed in terms of $\ww\eta1$, $\ww\eta2$ and $\ww\eta3$. To make the notation less cumbersome we
introduce the functions
\begin{equation*}
  \ww\teta i=\frac{\ww\eta i}{2^\sn{i+1}k^\sn2}=r^\sn{2i}+s^\sn{2i}+t^\sn{2i},\quad i=1,2,3\text{.}
\end{equation*}
Moreover we put
\begin{gather*}
  \al=\sqrt{\qq\teta13 - 3\ww\teta1\ww\teta2 + 2\ww\teta3},\quad\be=9\,\al^\sn2 - 5\,\qq\teta13 
  + 9\,\ww\teta1\,\ww\teta2\text{,}\\
  \psi=\be^\sn2 + 2\,{\left( \qq\teta12 - 3\,\ww\teta2 \right)}^3\quad\text{and}\\
  \ka=2^\sn{\frac13}\,{\left(\be+\sqrt\psi\right)}^\sn{\frac13} + 2\,\ww\teta1 -
  \frac{2^\sn{\frac23}\,\left(\qq\teta12-3\,\ww\teta2\right)} {{\left(\be + \sqrt\psi\right)}^{\frac13}}\text{.}
\end{gather*}
The potential is now given by
\begin{equation}
  \label{eq:Explicit-HK-pot-cohom-3-general-orbit}
  \rho=\frac{2\,\sqrt2\,k^\sn2}{\sqrt3}\,\left(\sqrt\ka+\sqrt{\frac{12\,\al}{\sqrt\ka}+6\,\ww\teta1-\ka}\right)\text{.}
\end{equation}
This was obtained by noticing that $\la=\frac{\rho}{4\,k^\sn2}=r+s+t$ is a solution of the equation
\begin{equation*}
 (\la^\sn2-\ww\teta1)^\sn2-\frac{8\,\al}{\sqrt 6}\,\la-2\left(\qq\teta12-\ww\teta2\right)=0\text{.}
\end{equation*}

\subsubsection{Relation with Hermitian Symmetric Spaces}
\label{section:Rel-w-Hermitian-Sym-Spaces}
Recall that we labelled the nilpotent orbits of cohomogeneity three according to whether or not the Lie span
$\List{X,X'}$ lies in a $\si$-invariant ${\slc}^\sn3$-subalgebra for an arbitrary element $X$. It will be shown in
\S~\ref{section:Cohom3-Orbit-in-SO7} that the orbit $(3,2^\sn2)$ in $\son7$ does not meet the requirement to be a
generic cohomogeneity-three orbit. As for the remaining orbits in Table~\ref{Table:All-cohom-three} we use an
interesting relation with Hermitian symmetric spaces to conf\mbox{}irm the generic cohomogeneity-three property.

\begin{remark}
  It is only the orbit in $\ww\LE7$ that calls for these considerations. In fact, the discussion in
  \S~\ref{sec:Examplex-of-Diagonalization} covers any cohomogeneity-three nilpotent orbit in a classical Lie algebra
  dif\mbox{}ferent from $\son7$.
\end{remark}

Let $\lie r\subset\liegc$ be a real semisimple Lie algebra with Cartan decomposition $\lie r=\lie k+\lp$ and
corresponding Cartan involution $\theta$. That is $\lie k=\lie r^\sn\theta$, $\lp=\lie r^\sn{-\theta}$ are such that
\begin{gather}
  \label{eq:CartanDecomp}
  \DL{\lie k\lie k}\subseteq\lie k,\quad\DL{\lie k\lp}\subseteq\lp,\quad\DL{\lp\lp}\subseteq\lie k
\end{gather}
and $\lie k+i\,\lp$ is a compact real form of $\lierc=\liegc$. In particular we choose $\lie r$ such that $\lie
k+i\,\lp=\lie g$ so that $\si$ becomes the conjugation map of $\lierc$ with respect to $\lie k+i\,\lp$. Let $\Lie
R\subset\Liegc$ be the connected Lie group with Lie algebra $\lie r$. Then $\LK$ is a Lie subgroup of $\LR$ with Lie
algebra $\lie k$ and, provided the center of $\LR$ is f\mbox{}inite, the homogeneous space $\LR/\LK$ is a
Riemannian symmetric space (of non-compact type).
\begin{definition}
  A pair $\left(\lie r,\ww H0\right)$ consisting of a real semisimple Lie algebra $\lie r$ with Cartan decomposition
  $\lie r=\lie k+\lp$ and an element $\ww H0$ in the center $Z(\lie k)$ of $\lie k$ such that the restriction of
  $\ww\ad{\ww H0}$ to $\lp$ is a complex structure, is called a \emph{semisimple Lie algebra of Hermitian type}.
\end{definition}
Notice, if such an element $\ww H0$ exists then $J=\ww\ad{\ww H0}$ clearly induces a $\LR$-invariant almost complex
structure on $\LR/\LK$. On the other hand, if $\LR/\LK$ is Hermitian it follows that $\lie r$ is of
Hermitian type (see eg.~\cite{Knapp,Sakate}).

For the rest of this section we assume $\left(\lie r,\ww H0\right)$ to be of Hermitian type; consequently $\liegc$ is no
longer random among the complex simple Lie algebras. With the assignment $\ww J0=\ww\ad{\ww H0}\ww|{\lp}$ the
complexif\mbox{}ied Lie algebra decomposes as $\lierc=\lie k^\sn\C+\lp^\sn++\lp^\sn-$, where $\lp^\sn\pm$ is
the $\pm i$ eigenspace of $\ww J0$ and $\lp^\sn-=\Overline{\lp^\sn+}$. From~\eqref{eq:CartanDecomp} it follows
that
\begin{gather*}
  \DL{\lie k^\sn\C\lp^\pm}=\lp^\pm,\quad\DL{\lp^\sn+\lp^\sn-}\subseteq\lie k^\sn\C,\quad\DL{\lp^\pm\lp^\pm}=0\text{.}
\end{gather*}
So $\lp^\pm$ are in fact abelian subalgebras of $\lierc$. Let $\lie h$ be a maximal abelian subalgebra of $\lie k$,
and notice that $\ww H0\in Z(\lie k)\subseteq\lie h$. Since $\ww J0Z\ne0$ for nonzero $Z\in\lp$, the centralizer
of $\ww H0$ in $\lie r$ equals $\lie k$. Thus $\lie h$ is a maximal abelian subalgebra of $\lie r$ and, because the
elements of $\liehc$ are semisimple, $\liehc$ is a Cartan subalgebra of $\lierc$. Let $\Phi$ denote the root system of
$\lierc$ relative to $\liehc$. Let moreover $\Phi^\sn c$ stand for the set of all compact roots, and $\Phi^\pms$ denote
positive/negative non-compact roots, respectively. In other words $\Phi=\Phi^{\sn c}\cup\Phi^\sn+\cup\Phi^\sn-$
with
\begin{equation*}
  \begin{split}
    \al(\ww H0)= \left\{\begin{array}{ll}
        0     &\    \al\in\Phi^{\sn c}\\
        \pm i &\    \al\in\Phi^\pm
      \end{array}\right.\text{.}
  \end{split}
\end{equation*}
We then have $\lp^\pms=\ww\sum{\al\in\Phi^\pms}\qq{\lie r}\al\C$ where $\qq{\lie r}\al\C$ is the root space of
$\lierc$ corresponding to $\al\in\Phi$, and one f\mbox{}inds that $\Overline Z$ lies in $\qq{\lie r}{-\al}\C$ whenever
$Z$ lies in $\qq{\lie r}\al\C$. Using that $\lp^\sn+$ is abelian it is easy to show that (see~\cite{Knapp})
\begin{lemma}
  Let $\lie r$ be a semisimple Lie algebra of hermitian type with Cartan decomposition $\lie r=\lie k+\lp$, and let
  $\lie h\subseteq\lie k$ be a Cartan subalgebra of $\lie r$. Then there exist $r\in\N$ and $r$ linearly independent
  positive non-compact roots $\ww\be1,\dots,\ww\be r$ relative to $\liehc$ such that $\ww\be j\pm\ww\be k$ is not a root
  for any $j\ne k$ and $\lie a=\qq\sum{j=1}r\C(\ww Ej+\ww{\Overline E}j)$
  is a maximal abelian subalgebra of $\lp$ ($\ww Ej$ is a nonzero $\ww\be j$-root vector).\qed
\end{lemma}
\begin{remark}
  The dimension of a maximal abelian subalgebra of $\lp$ is called the \emph{real rank} of $\lie r$.
\end{remark}

In the above Lemma we may in fact choose $\ww E j$ such that $\DL{\ww Ej\ww Ej\ww{\Overline E}j}=-2\ww Ej$. Seeing that
$\siEd{j}=-\ww{\Overline E}j$ we obtain a $\si$-invariant subalgebra isomorphic to $\slc$ via the complex span of
$\List{\ww Ej,\ww{\Overline E}j,\DL{\ww Ej\ww{\Overline E}j}}$. Thus $\lieac $ is a $\si$-invariant subalgebra
consisting of $r$ copies of $\slc$ where $r$ is the real rank of $\lie r$. It is a fact that any two maximal abelian
subalgebras of $\lp$ are $\Ad$-conjugate by an element $k\in\LK$. Hence
\begin{equation}
  \label{eq:p-as-conj-of-a}
  \lp=\ww\bigcup{k\in\LK} \ww\Ad k(\lie a)  
\end{equation}
and by the action of $\LK$ any element of $\lp^\sn+$ can be moved into a $\si$-invariant ${\slc}^\sn r$-subalgebra.

It turns out that the restricted root system of $\lie r$ (relative to $\lie a$) is of type $\ww{BC}r$ or $\ww Cr$ (see
eg.~\cite[p.~110]{Sakate}). The only system that arises in our context is $\ww Cr$ (see
Table~\ref{Table:Some-Hermitian-Sym-alg}); it is given by
\begin{equation}
  \label{eq:Root-sys-Cr}
  \List{\pm\ww\xi j\pm\ww\xi k\,,\,\pm2\,\ww\xi j\,|\,1\leqs j,k \leqs r\,,\,j\ne k}\text{,}
\end{equation}
where $\List{\ww\xi j\,|\,1\leqs j\leqs r}$ is the dual basis of $\List{\ww Xj=\ww Ej+\ww{\Overline E}j\,|\,1\leqs
  j\leqs r}$. Consider
the $j$'th $\slc$ subalgebra of $\lieac$ and decompose $\liegc$ under the action of this subalgebra. If this is a
highest root decomposition,
\begin{equation*}
 \liegc\cong\ww\slc j\op\qq{\lie k}j\C\op\left(\qq Sj1\ot\ww Vj\right)\text{.}
\end{equation*}
we see via~\eqref{eq:Root-sys-Cr} that
\begin{equation}
  \label{eq:Vj-in-highest-root-decomp}
  \ww Vj\cong\C^\nu\ot\left(\ww\bigoplus{k\ne q}\qq Sk1\right)\text{,}
\end{equation}
where $\nu\geqs1$ is the dimension of the root space corresponding to the restricted root $\pm\ww\xi j\pm\ww\xi k$.
\begin{table}[h]
  \centering
  \begin{tabular}{cc}
    \hline\hline
    \TS$\lie r\subset\liegc$ & $\lie k$\\
    \hline
    \TS$\su{3,3}\subset\sln6$ & $\R\op\su3\op\su3$\\
    \TS$\lie{sp}(3,\R)\subset\lie{sp}(3,\C)$ & $\R\op\su3$\\
    \TS$\lie{so}^\ast(12)\subset\son{12}$ & $\R\op\su6$\\ 
    \TS$\ww{\lie e}{7(-25)}\subset\qq{\lie e}7\C$ & $\R\op\ww{\lie e}6$\\
    \hline\hline
  \end{tabular}
    \caption[Selected simple Lie algebras of Hermitian type of real rank $3$]{Selected simple Lie algebras of Hermitian
      type of real rank $3$.}
    \label{Table:Some-Hermitian-Sym-alg}
\end{table}

We are now in position to verify Table~\ref{Table:All-cohom-three}. First, let $\lie r$ be one of the Hermitian
symmetric Lie algebras of Table~\ref{Table:Some-Hermitian-Sym-alg} and keep in mind the Beauville-data in
Tables~\ref{Table:Big-table-classical} \&~\ref{Table:Big-table-execptional} of the nilpotent orbit $\O$ of cohomogeneity
three in the corresponding complexif\mbox{}ied Lie algebra. It is then evident that the $\LK$-module $\lie n$ equals
$\lp^\sn+$ (or $\lp^\sn-$). Since any element $X$ of $\O$ can be moved into $\lie n$ by the action of $\G$ we conclude
via~\eqref{eq:p-as-conj-of-a} that $X$ can be put into a $\slT$-subalgebra. Thus $\O$ is a generic cohomogeneity three
orbit. In addition, since $X$ may be decomposed as in~\eqref{eq:decomp-X}, the isomorphism
\eqref{eq:Vj-in-highest-root-decomp} gives another argument for the existence of the one-parameter family of
hyperK\"ahler metrics with K\"ahler potentials. Second, let $\lie g$ be one of the Lie algebras $\lie{sp}(3+i)$,
$\su{6+i}$, $\So{12+i}$ with $i\geqs1$. A look at Tables~\ref{Table:Big-table-classical}
\&~\ref{Table:Big-table-execptional} reveals that, in spite of the growth of $\LK$ as $i$ increases, the action of $\LK$
on $\lie n$ is independent of $i\geqs0$. In other words, moving an element of $\lie n$ via the action of $\LK$ reduces
to the situation where $i=0$. But this is already accounted for by the appearance of a Lie algebra of Hermitian type. This
covers all the orbits in Table~\ref{Table:All-cohom-three}.

\subsection{The Special Cohomogeneity Three Orbit}
\label{section:Cohom3-Orbit-in-SO7}
Recall that $\son7$ possesses one nilpotent orbit $\O=(3,2^\sn2)$ of cohomogeneity three. Unfortunately this orbit
behaves dif\mbox{}ferently from the ones of \S~\ref{sec:general_cohom_three_orbits}. In fact for a typical element
$X\in\O$, $X$ and $X'$ span the whole ambient Lie algebra $\son7$ (see~\cite{MathematicaSO7-2}). Thus there is no
alternative left than computing the endomorphism $J$ on all of $\son7$. For this task we use \textsc{Mathematica} to
derive equations~\eqref{eq:J2-SO7-mathematica},~\eqref{eq:J2-SO7-a} and~\eqref{eq:J2-SO7-b} below. The complete
computations are collected in a notebook f\mbox{}ile~\cite{MathematicaSO7-2}. In principle, these equations are
computable by hand as well, using~\eqref{eq:J-cohom3} and~\eqref{eq:eta123-SO7-generic-X}.

\

Via the action of the compact group $G=\SO(7)$ it is possible to move any element $X$ of $\O$ into the subspace $\lie
n$. However, this can be improved. The weighted Dynkin diagram of $\O$ is $10\rtwoa1$ so there is a basis
$\List{\al,\be,\ga}$ of simple positive roots such that $\ww\ad H$ acts on $\qq{\lie g}\al\C$, $\qq{\lie g}\be\C$ and
$\qq{\lie g}\ga\C$ with eigenvalues $1$, $0$ and $1$ respectively. The eigenspaces of $\ww\ad H$ with eigenvalues larger
than one are $\liegc(2) = \qq{\lie g}{\al+\be+\ga}\C \op \qq{\lie g}{\be+2\ga}\C$ and $\liegc(3) = \qq{\lie
  g}{\al+\be+2\ga}\C \op \qq{\lie g}{\al+2\be+2\ga}\C$. The zero-eigenspace is $\liegc(0)=\ww\C+\op\ww\C-\op\slc$ so
$\LK=\ww{\Lie U(1)}+\ww{\LU(1)}-\SU(1)$ and one f\mbox{}inds that $\liegc(2)=\Lw+ + \Lm2$ and $\liegc(3)=\Lw+\Lw-\Se\Ct$
as a $\LK$-module. Thus, f\mbox{}ixing nonzero root vectors $\ww E\al\in\qq{\lie g}\al\C$ etc., the action of
$\LK\subset G$ allows us to write a generic element as $X=r\,\ww E{\al+\be+\ga} + s\,\ww E{\be+2\ga} + t\,\ww
E{\al+\be+2\ga}$ with $r,s,t>0$.

We choose to describe $\son7$ as the set of complex $(n\times n)$ matrices $Z$ such that $ZB+ZB^\sn t=0$ where
\begin{equation*}
  B=
  \begin{pmatrix}
    1&0&0\\
    0&0&\ww I3\\
    0&\ww I3&0
  \end{pmatrix}\text{.}
\end{equation*}
The diagonal elements form the Cartan subalgebra and the compact real form is $\si: Z \mapsto -{\Overline{Z^\sn t}}$. For
each positive root $\phi$ the root vector $\ww E\phi$ is chosen such that $\List{\ww H\phi,\ww E\phi,\ww F\phi}$ is a
standard basis of $\sln2$ with $\ww H\phi=\DL{\ww E\phi\ww F\phi}$, $\ww F\phi=-\siEd\phi$. Explicitly, $X$ takes
the matrix form
\begin{equation*}
  X=
  \begin{smatrix}
    0&0&0&0&\sqrt2r&0&0 \\
    -\sqrt2r&0&0&0&0&t&0& \\
    0&0&0&0&-t&0&s \\
    0&0&0&0&0&-s&0 \\
    0&0&0&0&0&0&0 \\
    0&0&0&0&0&0&0 \\
    0&0&0&0&0&0&0 \\
  \end{smatrix}\text{.}
\end{equation*}
The three invariant functions evaluated at $X$ are
\begin{equation}
 \label{eq:eta123-SO7-generic-X}
  \begin{split}
   \ww\eta1(X)& = 10\, \left( 2\,r^\sn2 + s^\sn2 + t^\sn2 \right),\quad
  \ww\eta2(X) = 20\, \left( 2\,r^\sn4 + 4\,r^\sn2\,t^\sn2 + {\left( s^\sn2 + t^\sn2 \right)}^2 \right)\\
   \text{and}&\quad \ww\eta3(X) = 40\, \left( 2\,r^\sn6 + 9\,r^\sn4\,t^\sn2 + 6\,r^\sn2\,t^\sn2 \left( s^\sn2 + t^\sn2
     \right) + {\left( s^\sn2 + t^\sn2 \right)}^\sn3 \right)\text{.}
  \end{split}
\end{equation}

Our aim is to use~\eqref{eq:J-cohom3} to compute the endomorphism $J^\sn2+1$ on a basis of the tangent space at
$X$ and then impose $J$ to be almost complex. The complex dimension of $\O$ is $12$ therefore the condition $J^\sn2+1=0$
implies $144$ new equations. A minority of these are nontrivial, and in particular one f\mbox{}inds that the K\"ahler
potential $\rho=\rho(r,s,t)$ must solve the following system of equations,
\begin{equation}
\label{eq:J2-SO7-mathematica}
  \begin{split}
  \ww\rho s \, \ww\rho{rs} +    \ww\rho t \, \ww\rho{rt} = \ww\rho s \, \ww\rho{st} +    \ww\rho t \, \ww\rho{tt} =
  \ww\rho r \, \ww\rho{rs} + 4\,\ww\rho t \, \ww\rho{st} = \ww\rho r \, \ww\rho{rt} + 4\,\ww\rho t \, \ww\rho{tt} =0\\
  \end{split}
\end{equation}
and
\begin{subequations}
  \begin{gather}
    \label{eq:J2-SO7-a}
    \begin{split}
      t \left( 100 - \qq\rho t2 \right) =\left( r\,\ww\rho r + s\,\ww\rho s \right) \ww\rho t
    \end{split}
    \\
    \label{eq:J2-SO7-b}
    \begin{split}
      t\left(r\,\ww\rho r - s\,\ww\rho s\right)=\left(4\,r^\sn2-s^\sn2\right)\ww\rho t\text{.}
    \end{split}
  \end{gather}
\end{subequations}
Consequently, $\qq\rho s2+\qq\rho t2$ is a function of $s$ and $\qq\rho r2+4\,\qq\rho
t2$ is a function or $r$. A combination of~\eqref{eq:J2-SO7-a} and~\eqref{eq:J2-SO7-b} gives an equation involving
quadratic terms of the f\mbox{}irst derivatives
\begin{equation}
  \label{eq:Combi-J2-a-b}
  s^\sn2\left(\qq\rho s2 + \qq\rho t2 - 100\right) = r^\sn2\left(\qq\rho r2 + 4\,\qq\rho t2-400\right)\text{.}
\end{equation}
Notice that the left hand side of~\eqref{eq:Combi-J2-a-b} is a function or $s$, whereas the right is a function of
$r$, so both must be constant. Thus
\begin{subequations}
  \begin{gather}
    \label{eq:SO7-rho2-st}
    \begin{split}
      \qq\rho s2+\qq\rho t2=100\left(1+\frac c{s^\sn2}\right)\text{,}
    \end{split}
    \\
    \label{eq:SO7-rho2-rt}
    \begin{split}
      \qq\rho r2+4\qq\rho t2=100\left(4+\frac c{r^\sn2}\right)
    \end{split}
  \end{gather}
\end{subequations}
for some constant $c$.

Eliminating $\ww\rho r$ and $\ww\rho s$ from the equations~\eqref{eq:J2-SO7-a},~\eqref{eq:J2-SO7-b}
and~\eqref{eq:SO7-rho2-st}, we obtain a quadratic equation in $\qq\rho t2$,
\begin{equation*}
    50\,t^\sn4=t^\sn2 \left(\left( 4\,r^\sn2 + s^\sn2 + t^\sn2 \right) + 2\,c \right)\qq\rho t2
    -\frac{{\left( 4\,r^\sn2 + s^\sn2 + t^\sn2 \right)}^\sn2 - 16\,r^\sn2 s^\sn2}{200}\qq\rho t4\text{.} 
\end{equation*}
In order for $\qq\rho t2$ to be real for all $r,s,t$ one needs that $c\geqs0$. Solving for $\ww\rho t$ leads to
\begin{equation}
  \label{eq:SO7-derivative-rho-t}
  \begin{split}
    \frac{50\,t^\sn2}{\qq\rho t2}=\vep\,\sqrt{4\,r^\sn2\,s^\sn2+c\,\left(4r^\sn2\,+s^\sn2+t^\sn2\right)
      +c^\sn2}+c +\frac{\left(4\,r^\sn2 + s^\sn2 + t^\sn2\right)}2
  \end{split}
\end{equation}
where $\chi=\pm1$ and $\vep=\pm1$.
Integrating, we f\mbox{}ind that $\O$ carries a one-parameter family of hyperK\"ahler metrics with K\"ahler potentials,
\begin{equation}
 \label{eq:SO7-Int-rho-t}
\begin{split}
  \rho(r,s,t;c)=&10\,\chi \,\Bigl( \sqrt{ 4\,r^2 + s^2 + t^2  +  2\, h(r,s,t) } + f(r,s) \\
  &- \sqrt c\, \log \left( h(r,s,t) + \sqrt c \, \sqrt{4\,r^2 + s^2 + t^2 + 2\, h(r,s,t)}\right)\Bigr)
\end{split}
\end{equation}
where
\begin{equation*}
  h(r,s,t)=c+\vep\,\sqrt{4\,r^2\,s^2 + \left( 4\,r^2 + s^2 + t^2 \right) \,c  + c^2}
\end{equation*}
and $f(r,s)$ is a function to be determined. The choice of the positive solution ($\chi=1$) forces one to take $\vep=1$,
since $\rho$ should be everywhere dif\mbox{}ferentiable.

To f\mbox{}ind $f$ one may write the right hand side of~\eqref{eq:SO7-Int-rho-t} as $\go I(r,s,t)+f(r,s)$ so that $\ww
fs=\ww\lim{t\rightarrow0}\left(\ww\rho s - \ww{\go I}s\right)$. The derivative $\ww\rho s$ can be found using the
equations~\eqref{eq:SO7-rho2-st} and~\eqref{eq:SO7-derivative-rho-t}. It turns out that
\begin{equation*}
  f(r,s)=\sqrt c\log\left(r\,s\right)
\end{equation*}
up to addition of some constant.

For $c=0$ the solution simplif\mbox{}ies somewhat,
\begin{equation}
  \label{eq:HK-pot-SO7-local-expression}
  \rho(r,s,t;0)=10\sqrt{\left(2\,r + s\right)^\sn2+t^\sn2}\text{.}
\end{equation}

Let us introduce a new invariant function on $\O$,
\begin{equation*}
  \ww\ze3(Z)=\frac15\,\ww\eta1(Z^\sn2)=\frac15\,\inp{Z^\sn2}{\si Z^\sn2}=\tr{ZZ\Overline{Z^\sn
  t}\Overline{Z^\sn t}},\quad Z\in\O
\end{equation*}
which is non-trivial since $\O$ consists of $3$--step nilpotent elements. Evaluated at $X$ we get that
$\ww\ze3(X)=4\,r^\sn4$. It is convenient to def\mbox{}ine $\ww\ze j(Z)=\ww\eta j(Z)/5$, $j=1,2$; notice that
$\ww\ze1(Z)=\tr{Z\Overline{Z^\sn t}}$ and $\ww\ze2(Z)=\tr{\DL{Z\Overline{Z^\sn
      t}}^\sn2}$. Now~\eqref{eq:HK-pot-SO7-local-expression} can be written as
\begin{equation}
  \label{eq:hyperKahler-pot-SO7-zero-constant}
  \rho=5\,\sqrt{2}\,\sqrt{\ww\ze1+2\sqrt{\ww\ze3}+2\sqrt{\qq\ze12 - \ww\ze2 - 2\,\ww\ze3}}\text{,}
\end{equation}
This is the unique hyperK\"ahler potential, because for $\la\in\C$ it is obvious that $\ww\ze1(\la
Z)=|\la|^\sn2\ww\ze1(Z)$ and $\ww\ze2(\la Z)=|\la|^\sn4\ww\ze2(Z)$, $i=2,3$ so $\rho(\la Z)=|\la|\rho(Z)$ (see
\cite{Brylinski}).

In fact it is possible to write down the one-parameter family of K\"ahler potentials in terms of the globally
def\mbox{}ined $G$-invariant functions $\ww\ze1$, $\ww\ze2$ and $\ww\ze3$,
\begin{equation}
  \label{eq:SO7-one-parameter-familiy-Kahler-pot}
     \rho=\frac{10}{\sqrt2} \left(\sqrt{\ww\ze1 +2\,\sqrt{\ww\ze3} + 4\,h}-\sqrt{2\,c}\,\log\left( \frac{h +
      \frac c{\sqrt2}\,\sqrt{\ww\ze1 +2\,\sqrt{\ww\ze3}+4\,h}}{\sqrt{\qq\ze12-\ww\ze2-2\,\ww\ze3}}\right)\right)\text{,}
\end{equation}
where now
\begin{equation*}
  h = c + \sqrt{\frac14\left( \qq\ze12 - \ww\ze2 - 2\,\ww\ze3\right) + 
    c\,\left(\frac12\ww\ze1 + \sqrt{\ww\ze3}\right) + c^\sn2}\text{.} 
\end{equation*}

\begin{remark}
  \label{remark:HK-pot-SO7-orbit-322}
  Using hyperK\"ahler quotients Kobak \& Swann~\cite{KobakSwann:HK-pot-via-finite-dim-quotients} found the invariant
  hyperK\"ahler potential for the hyperK\"ahler structure on the nilpotent orbit $(\O_{(3,2^\sn2)},\Si)\subset\son n$.
  The expression given in ~\cite{KobakSwann:HK-pot-via-finite-dim-quotients} matches what we found
  in~\eqref{eq:hyperKahler-pot-SO7-zero-constant}.
\end{remark}

\begin{remark}
  Kobak \& Swann showed in \cite{KobakSwann:Classical-Nilpotent-Orbits-HKQ} the existence of a one-to-one correspondence
  between the nilpotent orbits $\O_{(2^\sn4)}\subset\son8$ and $\O_{(3,2^\sn2)}\subset\son7$. The nilpotent orbit
  $\O_{(2^\sn4)}\subset\son8$ is of cohomogeneity two and admits a one-parameter family of $\SO(8)$-invariant
  hyperK\"ahler metrics with $\SO(8)$-invariant K\"ahler potentials which includes the unique hyperK\"ahler potential
  (see \cite{KobakSwann:HK-pot-in-cohom-two}). Reducing the symmetry group to $\SO(7)$ the above shows that there are no
  extra solutions, even though the cohomogeneity changes to three.
\end{remark}

\subsection{The Main Theorem}
In \S~\ref{sec:general_cohom_three_orbits} and \S~\ref{section:Cohom3-Orbit-in-SO7} we considered hyperK\"ahler metrics
with K\"ahler potential on all nilpotent orbits of cohomogeneity three under the action of the compact group $\G$. We
proved that the majority of orbits allowed only one such metric. Exceptions were found either when the Beauville bundle
is a vector bundle over one of the following (compact) Hermitian symmetric spaces
\begin{equation*}
  \frac{\SU(6)}{\LS(\ww\LU3\times\ww\LU3)},\quad\frac{\SP(3)}{\LU(3)},\quad\frac{\SO(12)}{\LU(6)},
  \quad\frac{\ww\LE7}{\ww\LE6\LU(1)}
\end{equation*}
or in the special case $\G=\SO(7)$. To summarize,

\begin{theorem}
  \label{Theorem:main-claim}
  Suppose $\G$ is a compact simple Lie group and $\O$ is a nilpotent orbit in $\liegc$ of cohomogeneity three. Then
  $\O$, endowed with the Kirillov-Kostant-Souriau complex symplectic form $\Si$, admits a unique $\G$-invariant
  hyperK\"ahler metric with $\G$-invariant hyperK\"ahler potential. In fact, the metric is the unique $\G$-invariant
  hyperK\"ahler metric of $(\O,\Si)$ with a K\"ahler potential unless $\lieg$ is one of $\lie{sp}(3)$, $\su6$,
  $\So{12}$, $\ww{\lie e}7$, $\So7$, in which case the metric lies in a one-parameter family of hyperK\"ahler metrics
  with K\"ahler potential.
  
  The unique hyperK\"ahler potential is given by~\eqref{eq:Explicit-HK-pot-cohom-3-general-orbit} for the generic
  cohomogeneity-three orbits, and by~\eqref{eq:SO7-one-parameter-familiy-Kahler-pot} for the special orbit.
\end{theorem}

\begin{proof}
  The only thing left to consider is the issue of the $\G$-invariant hyperK\"ahler potential. First, consider the
  generic cohomogeneity-three orbits. By Swann~\cite{Swann:HK-quaternionicKahler} it is known that any hyperK\"ahler
  metric on $\O$ admits a hyperK\"ahler potential. Consequently, when the $\G$-invariant K\"ahler potential $\rho$ is
  unique $(c=0)$ there is nothing to prove. Consider now the cases where $\O$ admits a one-parameter family of
  $\G$-invariant K\"ahler potentials. Now any $X\in\O$ lies in a real $\slT$ subalgebra, so if $\rho$ is a hyperK\"ahler
  potential it must restrict to a hyperK\"ahler potential of the regular orbit $\ww\O{reg}$ of $\slT$. But the
  hyperK\"ahler metric on $\ww\O{reg}$ is a product of three hyperK\"ahler metrics on the minimal orbit of $\slc$, and
  on each factor the hyperK\"ahler potential is unique ($c=0$), see Remark~\ref{Remark:unique_hKpot}.
  At last, for the special cohomogeneity three orbit we showed that the hyperK\"ahler potential is given
  by~\eqref{eq:hyperKahler-pot-SO7-zero-constant}.
\end{proof}

\begin{acknowledgements} This work is part of the author's Ph.D.~thesis at the University of Southern Denmark. He wishes
  to thank A.~F.~Swann for support and encouragement. The author is a member of \textsc{Edge} Research Training Network
  \textsc{hprn-ct-\oldstylenums{2000}-\oldstylenums{00101}}, supported by the European Human Potential Programme.
\end{acknowledgements}

\providecommand{\bysame}{\leavevmode\hbox to3em{\hrulefill}\thinspace}


\begin{thebibliography}{10}

\bibitem{Beauville:fano}
A.~Beauville, \emph{Fano contact manifolds and nilpotent orbits}, Comment.
  Math. Helv. \textbf{73} (1998), 566--583.

\bibitem{Besse:Einstein}
A.~L. Besse, \emph{Einstein {Manifolds}}, Ergebnisse der Mathematik und ihrer
  Grenzgebiete, 3. Folge, vol.~10, Springer-Verlag, Berlin Heidelberg, 1987.

\bibitem{Bourbaki1982}
N.~Bourbaki, \emph{Groupes et algébres de {Lie}, ch. 9}, Masson, Paris, 1982.

\bibitem{Bredon}
G.~E. Bredon, \emph{Introduction to {Compact} {Transformation} {Groups}}, Pure
  and applied mathematics, no.~46, Academic Press, New York, 1972.

\bibitem{Brylinski}
R.~Brylinski, \emph{{Instantons and Kaehler Geometry of Nilpotent Orbits}},
  "Representation Theories and Algebraic Geometry", Ed. A. Broer, NATO ASI
  Series, Kluwer, 1998, pp.~85--125.

\bibitem{Calabi:Kaehler}
E.~Calabi, \emph{Métriques kähleriennes et fibrés holomorphes}, Ann. Scient.
  {\'E}c. Norm. Sup. \textbf{12} (1978), 269--294.

\bibitem{CollingwoodMcGovern}
D.~H. Collingwood and W.~M. McGovern, \emph{Nilpotent {Orbits} in {Semisimple}
  {Lie} {Algebras}}, Van Nostrand Reihnhold, New York, 1993.

\bibitem{DancerSwann:HK-metrics-of-cohom-one}
A.~S. Dancer and A.~F. Swann, \emph{Hyperk\"ahler metrics of cohomogeneity
  one}, J. Geom. and Phys. \textbf{21} (1997), 218--230.

\bibitem{DancerSwann:qK-cohom-one}
\bysame, \emph{Quaternionic {K\"ahler} manifolds of cohomogeneity one}, Int. J.
  Math. \textbf{10} (1999), no.~5, 541--570.

\bibitem{GriffithsHarris}
P.~Grif\mbox{}f\mbox{}iths and J.~Harris, \emph{Principles of {Algebraic} {Geometry}}, John
  Wiley and Sons, New York, 1978.

\bibitem{Hitchin:monopoles}
N.~J. Hitchin, \emph{Monopoles, minimal surfaces and algebraic curves},
  Séminaire de Mathématiques Supérieures, 105, Presses de l'{U}niversit{\'e} de
  Montr{\'e}al, Montreal, PQ, 1987.

\bibitem{Knapp}
A.~W. Knapp, \emph{Lie {Groups} {Beyond} an {Introduction}}, Progress in
  Mathematics, no. 140, Bitkh\"auser, Boston, Mass., 1996.

\bibitem{KobakSwann:Classical-Nilpotent-Orbits-HKQ}
P.~Z. Kobak and A.~F. Swann, \emph{Classical nilpotent orbits as
  hyper{K\"ahler} quotients}, Int. J. Math. \textbf{7} (1996), no.~2, 193--210.

\bibitem{KobakSwann:HK-geo-wolf-spaces}
\bysame, \emph{The {HyperK\"ahler} {Geometry} {Associated} to {Wolf} {Spaces}},
  Bollettino U. M. I. \textbf{\textnormal{(8)} 4-B} (2001), 587--595.

\bibitem{KobakSwann:HK-pot-in-cohom-two}
\bysame, \emph{Hyper{K\"ahler} potentials in cohomogenety two}, J. reine angew.
  Math. \textbf{531} (2001), 121--139.

\bibitem{KobakSwann:HK-pot-via-finite-dim-quotients}
\bysame, \emph{{HyperK\"ahler} {Potentials} via {Finite-Dimensional}
  {Quotients}}, Geom. Dedicata \textbf{88} (2001), 1--19.

\bibitem{Kronheimer}
P.~B. Kronheimer, \emph{Instantons and the geometry of the nilpotent variety},
  J. Dif\mbox{}ferential Geom. \textbf{32} (1990), 473--490.

\bibitem{Sakate}
I.~Satake, \emph{{Algebraic} {Structures} of {Symmetric} {domains}},
  Publications of the Mathematical Society of Japan, vol.~14, Princeton
  University Press, 1980.

\bibitem{Swann:HK-quaternionicKahler}
A.~F. Swann, \emph{Hyperk\"ahler and quaternionic {K\"ahler} geometry}, Math.
  Ann \textbf{289} (1991), 421--450.

\bibitem{MathematicaSO7-2}
M.~Villumsen, \emph{The special cohomogeneity three orbit},
  \url{http://www.imada.sdu.dk/~mdv/so7/}, 2004.

\end{thebibliography}
\end{document}